\pdfoutput=1
\RequirePackage{ifpdf}
\ifpdf 
\documentclass[pdftex]{sigma}
\else
\documentclass{sigma}
\fi

\numberwithin{equation}{section}

\newtheorem{Theorem}{Theorem}[section]
\newtheorem{Lemma}[Theorem]{Lemma}
{ \theoremstyle{definition}
\newtheorem{Definition}[Theorem]{Definition}
\newtheorem{Remark}[Theorem]{Remark} }

\begin{document}

\allowdisplaybreaks

\renewcommand{\thefootnote}{$\star$}

\newcommand{\arXivNumber}{1603.07220}

\renewcommand{\PaperNumber}{076}

\FirstPageHeading

\ShortArticleName{$(D+1)$-Colored Graphs -- a Review of Sundry Properties}

\ArticleName{$\boldsymbol{(D+1)}$-Colored Graphs -- a Review\\ of Sundry Properties\footnote{This paper is a~contribution to the Special Issue on Tensor Models, Formalism and Applications. The full collection is available at \href{http://www.emis.de/journals/SIGMA/Tensor_Models.html}{http://www.emis.de/journals/SIGMA/Tensor\_{}Models.html}}}

\Author{James P.~RYAN}

\AuthorNameForHeading{J.P.~Ryan}

\Address{Institute for Mathematics, Astrophysics and Particle Physics, Radboud University,\\ Nijmegen, The Netherlands}
\Email{\href{mailto:j.ryan@science.ru.nl}{j.ryan@science.ru.nl}}

\ArticleDates{Received March 18, 2016, in f\/inal form July 25, 2016; Published online August 02, 2016}

\Abstract{We review the combinatorial, topological, algebraic and metric properties supported by $(D+1)$-colored graphs, with a focus on those that are pertinent to the study of tensor model theories. We show how to extract a limiting continuum metric space from this set of graphs and detail properties of this limit through the calculation of exponents at criticality.}


\Classification{60C05; 81Q35; 81T17; 83C27}

\renewcommand{\thefootnote}{\arabic{footnote}}
\setcounter{footnote}{0}

\section{Introduction}

This review focusses on the set of $(D+1)$-colored graphs, examining their combinatorial, topological, algebraic and metric properties. Such graphs form a subset of $(D+1)$-regular bipartite multigraphs distinguished by their admission of a specif\/ic edge-labelling. These labels or \emph{colors} are the key to their rich structure and the subsequent structural analysis.

Colored graphs arise mainly in two areas of the literature. Historically, they were f\/irst deve\-loped as a graph-theoretic tool to address challenges in piecewise-linear topology. See \cite{ferri-review} for an interesting review. More recently, they have experienced a resurgence of attention in the context of quantum gravity. As the Feynman graphs of tensor models/tensorial group f\/ield theories \cite{razvan-review, daniele-review}, they may form the groundwork for the systematic def\/inition of a non-trivial, physically interesting probability measure over quantum geometries.

From the tensor model perspective, the reason to concentrate a review purely on graph-theoretic properties is clearly motivated. Despite their striking similarities, tensors lack many of the powerful tools that are available for matrices \cite{matrix}. Thus, tensor model calculations rely for the moment quite heavily on one's ability to control and analyze the properties of suf\/f\/iciently interesting subsets of $(D+1)$-colored graphs.

The review falls into two parts. Sections \ref{sec:top} and \ref{sec:alg} introduce and analyse the combinatorial, topological and algebraic properties on the whole set of $(D+1)$-colored graphs. Afterwards, Section \ref{sec:mel} narrows the scope to the melonic subset, which permits a more in-depth analysis of their metric properties.

In Section \ref{sec:top}, we shall begin by carefully detailing the structure of $(D+1)$-colored graphs. From there, we use the colors to construct a topology on each graph. As a by-product of this construction, we demonstrate that such graphs encode $D$-dimensional simplicial pseudomanifolds. Thereafter, we show that the set of $(D+1)$-colored graphs admits several equivalence structures, one of which is known as combinatorial core equivalence. Each such class has preferred representatives, known as combinatorial core graphs and we present an algorithm that reduces a given graph to an equivalent core graph. To f\/inish this section, we identify a graph's jackets: specif\/ic embedded Riemann surfaces. They can be used to def\/ine an important combinatorial invariant, known as the degree. We detail some properties of the degree and its signif\/icance for the identif\/ication of melonic spheres.

In Section \ref{sec:alg}, we def\/ine the Lie algebraic structure supported by $(D+1)$-colored graphs. In Section \ref{sec:mel}, we focus entirely on one combinatorial core equivalence class: the (rooted) melonic graphs. We examine a metric structure thereon and most interestingly a continuum limit within this set of graphs. This limit has associated numbers, known as exponents, that help identify the limiting continuum metric space. We focus on three exponents: the susceptibility, Hausdorf\/f dimension and spectral dimension. All three indicate that continuum metric space coincides with that of branched polymers spacetimes.

As hinted above, this review rarely mentions its roots in either piecewise linear topology or tensor models. However, we take more care to highlight connections to the latter, while rarely, if ever, mention the former. Having said that, this is by no means a comprehensive review of all tensor model inspired, graph-theoretic properties. We make no ef\/fort to detail the properties of $(D+1)$-colored graphs with additional labels, which may be identif\/ied as matter \cite{ising-potts+,ising-potts}, dual-weighting \cite{dual-weighting}, richly-geometric \cite{group-field,group-field+}. We refrain from mentioning graph-theoretic properties pertaining to other core equivalence classes \cite{double-scaling+,double-scaling} or to multi-orientiable models \cite{multi-orient+,multi-orient}.

\subsection{Notation}

With such intricately labelled structures, a careful notation is required to decode the richly layered information. The colors associated to the edges are drawn from $\{0,1,\dots,D\}$. We shall be often interested in subsets $\{i,j,\dots,k\}$ and their complements $\{\widehat i,\widehat j,\dots, \widehat k\}= \{0,\dots, D \} \setminus \{ i,j, \dots, k \}$.

\section{Combinatorial and topological properties}\label{sec:top}

\subsection{Closed, open and boundary graphs}

Colored graphs are regular bipartite graphs that admit a specif\/ic labelling of the edges. Such graphs come in three varieties -- closed, open and boundary. Such qualif\/iers hint at a topological structure that we shall detail in due course. For the moment, we shall content ourselves with providing the following rather dry def\/initions:
\begin{Definition}[closed graph]\label{def:closed-graph}
A {closed $(D+1)$-colored graph} is a graph ${\mathcal G} = ({\mathcal V},{\mathcal E})$ with vertex set ${\mathcal V}$ and edge set ${\mathcal E}$ such that:
\begin{itemize}\itemsep=0pt
\item ${\mathcal V}$ is bipartite, that is, there is a partition of the vertex set ${\mathcal V} = V \cup \bar V$, such that for any element $l\in{\mathcal E}$, then $l = \{v,\bar v\}$ where $v\in V$ and $\bar v\in\bar V$. Their cardinalities satisfy $|{\mathcal V}| = 2|V| = 2|\bar V|$.
\item The edge set is partitioned into $D+1$ subsets ${\mathcal E} = \bigcup_{i =0}^{D} E^i$, where $E^i$ is the subset of edges with color $i$.
\item It is $(D+1)$-regular (i.e., all vertices are $(D+1)$-valent) with all edges incident to a given vertex having distinct colors.
\end{itemize}
\end{Definition}

\begin{figure}[t] \centering
 \includegraphics[scale = 1]{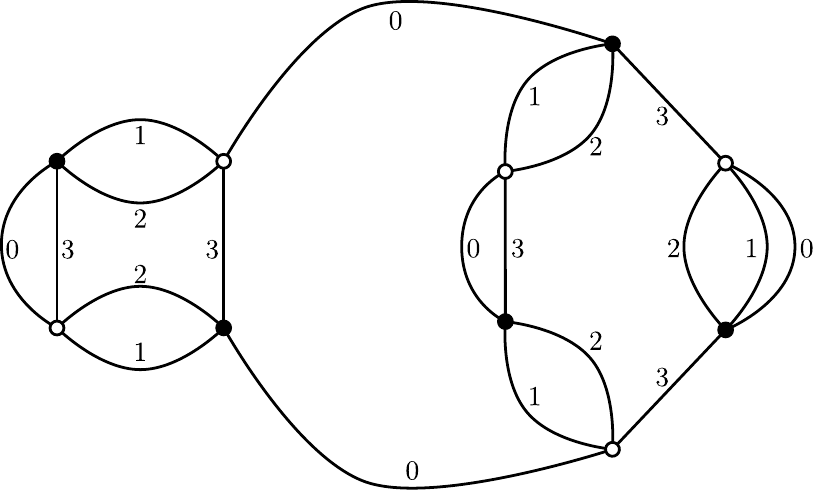}
 \caption{A closed $(D+1)$-colored graph ($D = 3$).}\label{fig:closed-graph}
\end{figure}

The elements $v\in V$ ($\bar v \in \bar{V}$) are commonly referred to as the positive (negative) vertices. Given a projection of the graph onto a~plane, one arranges the colors consistently in a clockwise (anti-clockwise) manner. Moreover, it is worth noticing that the bipartition induces an orientation on the edges, say from~$v$ to~$\bar v$. See Fig.~\ref{fig:closed-graph} for an example.

\begin{Definition}[open graph] \label{def:open-graph}
 An {open $(D+1)$-colored graph} is a graph ${\mathcal G}$ satisfying some additional constraints:
\begin{itemize}\itemsep=0pt
\item It is bipartite, that is, there is a partition of the vertex set ${\mathcal V} = V \cup \bar V$, such that for any element $l\in{\mathcal E}$, then $ l = \{v,\bar v\}$ where $v\in V$ and $\bar v\in\bar V$. Moreover, their cardinalities also satisfy $|{\mathcal V}| = 2|V| = 2|\bar V|$.

\item The positive vertices are of two types $V = V_{\rm int} + V_\partial$, where $V_{\rm int}$ is the set of $(D+1)$-valent {\it internal} vertices and the elements of $V_\partial$ are 1-valent {\it boundary} vertices. A similar distinction holds for negative vertices.

\item The edge set is partitioned into $D+1$ subsets ${\mathcal E} = \bigcup_{i=0}^D E^i$, where $E^i$ is the subset of edges with color $i$. Furthermore, each $E^i = E^i_{\rm int} \cup E^i_{\rm ext}$, such that {\it internal} edges $E^i_{\rm int}$ join two internal vertices, while {\it external} edges $E^i_{\rm ext}$ join an internal vertex to a boundary vertex.
\end{itemize}
\end{Definition}
\begin{Remark}In this text, we consider only \emph{connected} closed and open graphs. This is not just for the convenience of subsequent analysis, but is well motivated by their usage in tensor model theories; the terms in the perturbative expansion of tensor model cumulants are neatly labelled by connected graphs.
\end{Remark}

\begin{figure}[t]\centering
\includegraphics[width=15cm]{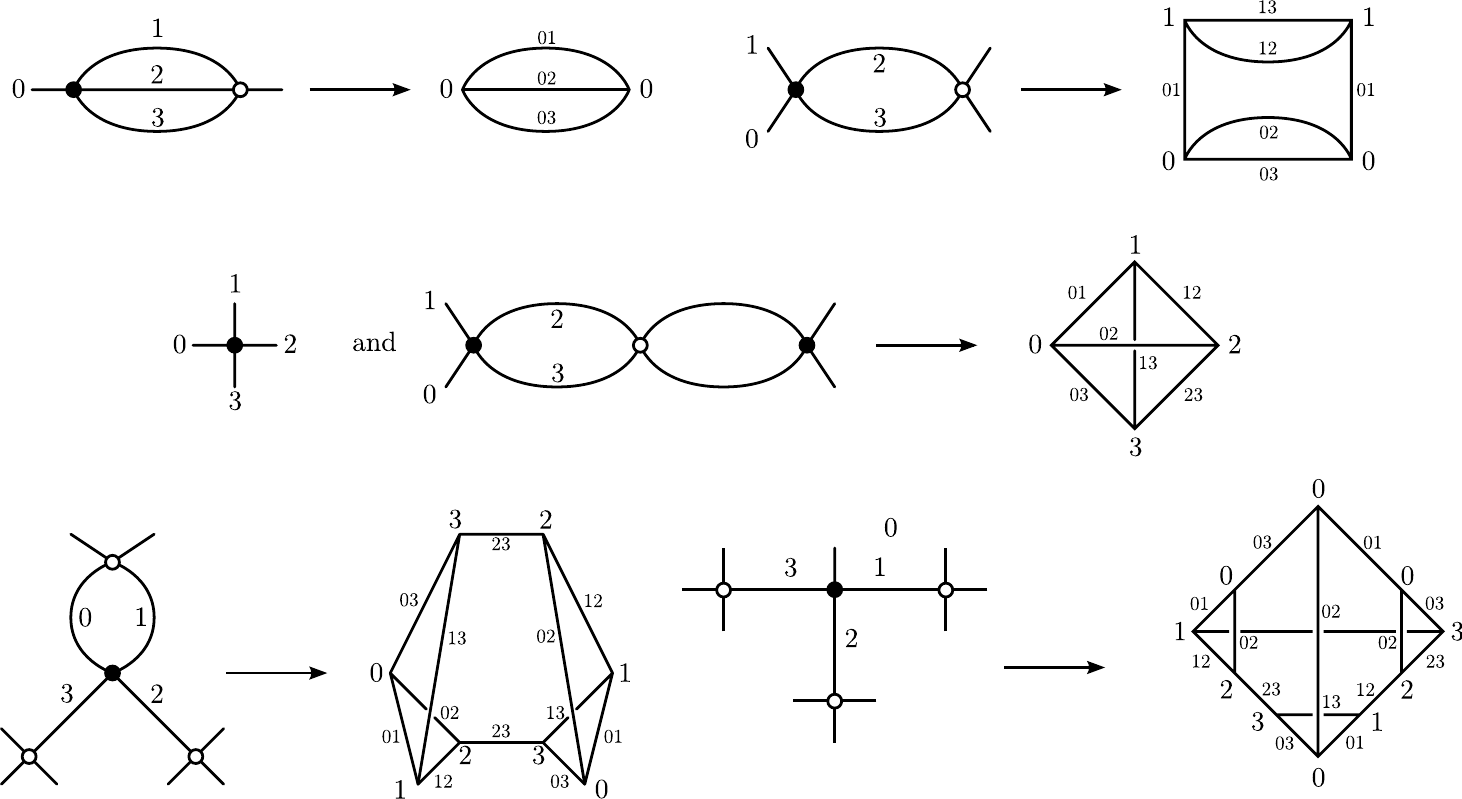}
\caption{Open 4-color graphs and their corresponding boundary graphs.}\label{fig:boundgr}
\end{figure}

The open graphs induce a boundary graph structure \cite{PolyColor} as follows:

\begin{Definition}[boundary graph]\label{def:boundary-graph}
 The {boundary graph ${\mathcal G}_{\partial}$} of an open $(D+1)$-colored graph~${\mathcal G}$ comprises of:
 \begin{itemize}\itemsep=0pt
 \item the vertex set ${\mathcal V}_\partial = V_{\partial}\cup \bar V_{\partial}$. We stress that it is not bipartite with respect to this splitting. The vertices inherit the color from the external edges of ${\mathcal G}$ upon which they lie, so that a~more appropriate partition is ${\mathcal V}_\partial = \bigcup_{i=0}^D {\mathcal V}_\partial^i$, where ${\mathcal V}_\partial^i$ denotes the set of boundary vertices with color $i$.
 \item the edge set ${\mathcal E}_\partial = \mathop{\bigcup_{i\neq j}} E^{ij}_\partial$, where $l^{ij} = \{v, w\} \in E^{ij}_\partial$ exists if there is a bi-colored path from $v$ to $w$ in ${\mathcal G}$ consisting of colors~$i$ and~$j$. Thus, the lines $E^{ij}_\partial$ inherit the colors of the path in ${\mathcal G}$.
 \end{itemize}
\end{Definition}
A cursory investigation of these boundary graphs reveals that they possess a number of additional properties. The edge $l^{ij} \in E_\partial^{ij}$ can only exist if $ l^{ij} = \{v^i,w^i\}$ or $\{v^i,w^j\}$ or $\{v^j,w^i\}$ or $\{v^j,w^j\}$, where $v^i, w^i \in {\mathcal V}_\partial^i$ and $v^j, w^j\in{\mathcal V}_\partial^j$. Each boundary vertex is $D$-valent and for $v^i\in {\mathcal V}_\partial^i$, the incident boundary edges are $l^{ij}$ where $j = \hat i$. Several examples are presented in Fig.~\ref{fig:boundgr}.

\begin{Remark}A priori, the boundary graph ${\mathcal G}_{\partial}$ is a very dif\/ferent beast from the initial graph~${\mathcal G}$; afterall, it has colored vertices and bicolored edges. Moreover, a connected open graph can have a boundary with several connected components. However, from a topological perspective, these disparate structures dovetail elegantly. We shall return to this later.
\end{Remark}

\subsection{Cellular structure and pseudomanifolds}

Such heavily labelled graphs display a beautiful hierarchy that is key to their utility \cite{lost, color}. Importantly, it supports a $D$-dimensional topological structure. We shall now expose this for a~generic $(D+1)$-colored graph~${\mathcal G}$.
\begin{Definition}[$d$-bubble]\label{def:bubble} The $d$-bubble of ${\mathcal G}$ is the maximally connected subgraph comprising of edges with~$d$ f\/ixed colors.
\end{Definition}

Obviously, $d\in\{0,\dots, D\}$. The $d$-bubbles are denoted by ${\mathcal B}^{i_1\dots i_d}_{ (\rho) }$; the color indices are ordered $i_1<i_2<\dots <i_d$ to uniquely identify the particular \emph{species} of $d$-bubble, while $\rho$ distinguishes connected components of the same species. We denote the number of $d$-bubbles of the graph by~${\mathcal B}^{[d]}$. One should note that the 0-bubbles are the vertices of~${\mathcal G}$, the 1-bubbles are the edges of~${\mathcal G}$. The $2$-bubbles are the faces of~${\mathcal G}$.

\begin{figure}[t]\centering
 \includegraphics[scale = 0.7]{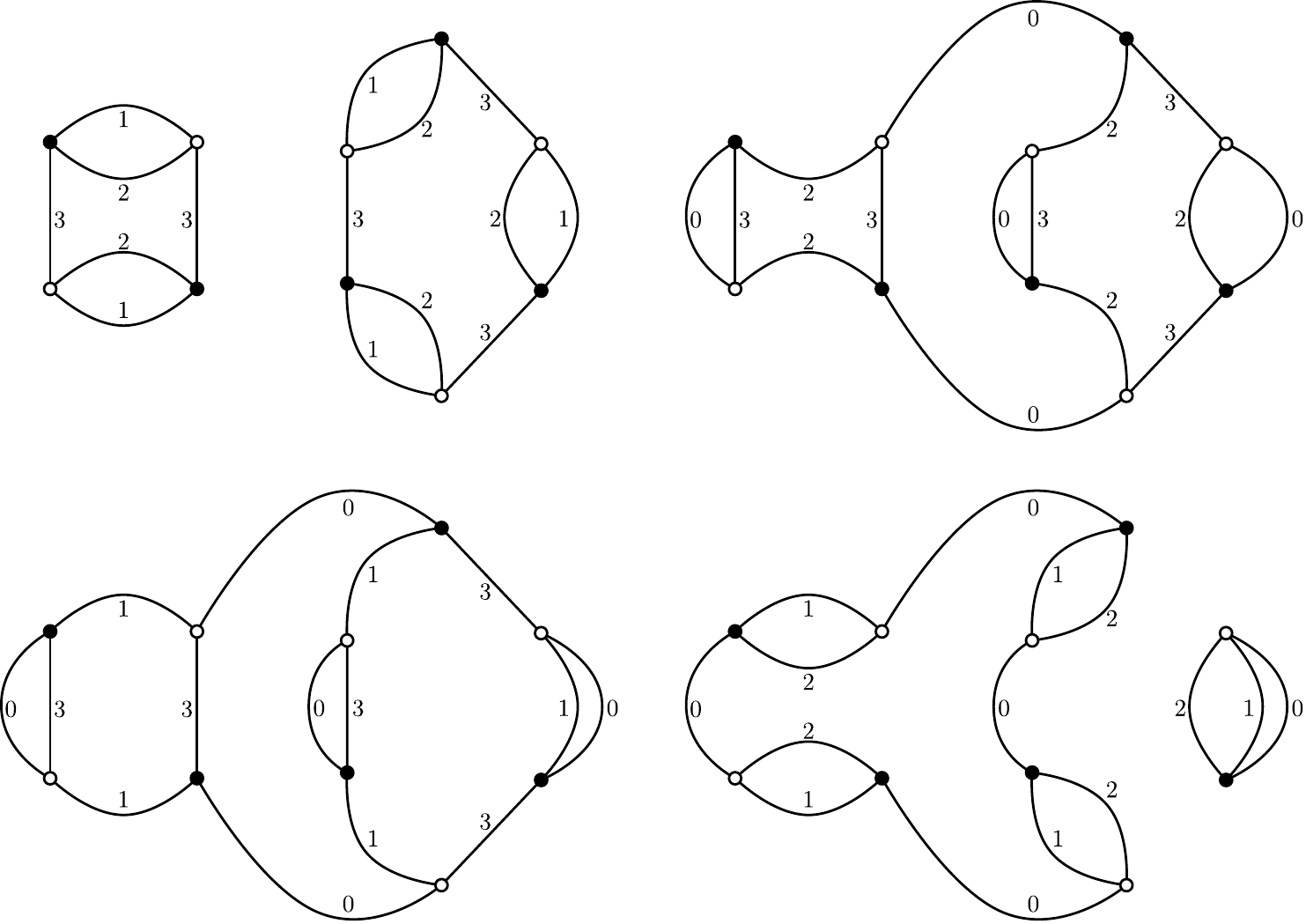}
\caption{$3$-bubbles of a graph in $D=3$.}\label{fig:exemplubub}
\end{figure}

In Fig.~\ref{fig:exemplubub}, we present the $3$-bubbles associated to the graph given in Fig.~\ref{fig:closed-graph}. The $3$-bubbles are indexed by the colors of their lines, namely from left to right $123$, $023$, $013$ and $012$. The $2$-bubbles are the subgraphs catalogued by color pairs $01$, $02$, $03$, $12$, $13$ and $23$. The $1$-bubbles are the lines $0$, $1$, $2$ and $3$, while the $0$-bubbles are the vertices.

\textbf{Constructing the dual complex.} These $d$-bubbles are key to def\/ining a topology, in particular, a $D$-dimensional cellular complex structure. To this end, we shall construct the dual \emph{finite abstract simplicial complex}. Quite clearly, the graph complex and the dual complex are the same topological space. Note that we write ${\mathcal H} \subset {\mathcal G}$, if ${\mathcal H}$ is a subgraph of~${\mathcal G}$. To construct the dual complex~\cite{lost}:
\begin{itemize}\itemsep=0pt
 \item We f\/irst assemble all the $D$-bubbles of ${\mathcal G}$ into a set~$A$:
\begin{gather*}
 A = \Bigl\{ {\mathcal B}^{\widehat i}_{(\rho)}\colon i \in \{ 0,\dots, D \},\; \rho \;\textrm{unrestricted} \Bigr\} .
 \end{gather*}
 \item For each $(D+1-d)$-bubble within ${\mathcal G}$, $d\in \{1,\dots, D+1\}$, we form the following subset of~$A$:
\begin{gather*}
 \sigma^{ {\mathcal B}^{ \widehat i_1 \widehat i_2\dots \widehat i_d}_{ (\kappa) } } =
 \Bigl\{ {\mathcal B}^{\widehat i_k}_{(\rho)} \colon {\mathcal B}^{ \widehat i_1 \widehat i_2\dots \widehat i_d}_{ (\kappa) }
 \subset {\mathcal B}^{\widehat i_k}_{(\rho)} , \; k\in \{1,\dots, d\} \Bigr\} \subset A .
 \end{gather*}
In fact, for a given $(D+1-d)$-bubble ${\mathcal B}^{ \widehat i_1 \widehat i_2\dots \widehat i_d}_{ (\kappa) }$, there is, for each $k\in\{1,\dots,d\}$, a unique $D$-bubble ${\mathcal B}^{\widehat i_k}_{(\rho)} \supset {\mathcal B}^{ \widehat i_1 \widehat i_2\dots \widehat i_d}_{ (\kappa) }$. This is the maximal connected component (in ${\mathcal G}$) obtained by starting from ${\mathcal B}^{ \widehat i_1 \widehat i_2\dots \widehat i_d}_{ (\kappa) }$ and adding edges of all colors except~$i_k$. Thus, the cardinality of~$\sigma^{ {\mathcal B}^{ \widehat i_1 \widehat i_2\dots \widehat i_d}_{ (\kappa) } }$ is~$d$.
\item As a result, any subset $\tau \in \sigma^{ {\mathcal B}^{ \widehat i_1 \widehat i_2\dots \widehat i_d}_{ (\kappa) } } $ is indexed by a choice of subset $S \subset \{ 1,\dots, d \} $:
 \begin{gather*}
 \tau = \Bigr\{ {\mathcal B}^{\widehat i_k}_{(\rho)} \,|
 \, {\mathcal B}^{ \widehat i_1 \widehat i_2\dots \widehat i_d}_{ (\kappa) }
 \subset {\mathcal B}^{\widehat i_k}_{(\rho)} ,\; k \in \{ 1,\dots, d \} \setminus S \Bigl\} .
 \end{gather*}
So $\tau = \sigma^{{\mathcal B}^{ \widehat i_1 \widehat i_2\dots \widehat i_{d-|S|}}_{ (\xi) }}$, where ${\mathcal B}^{ \widehat i_1 \widehat i_2\dots \widehat i_{d-|S|}}_{ (\xi) }$ is the unique subgraph obtained by adding the colors $i_s$, $s\in S$, to the subgraph ${\mathcal B}^{ \widehat i_1 \widehat i_2\dots \widehat i_d}_{ (\kappa) }$.

\item The f\/inal crucial detail is that the sets $\sigma^{{\mathcal B}^{ \widehat i_1 \widehat i_2\dots \widehat i_{d}}_{ (\kappa) }}$ are the $(d-1)$-simplices of a f\/inite abstract simplicial complex, def\/ined as
\begin{gather*}
\Delta = \Big\{ \sigma^{ {\mathcal B}^{\hat i_1\hat i_2\dots \hat i_d}_{ (\kappa) } } \, | \,
{\mathcal B}^{\hat i_1\hat i_2\dots \hat i_d}_{ (\kappa) } \subset {\mathcal G}\Big\}.
 \end{gather*}
It is straightforward to verify the def\/ining property of an abstract simplicial complex: that for all $\sigma \in \Delta$ and $\tau \subset \sigma$, then $\tau \in \Delta$.
The cardinality of $\sigma \in \Delta$ is $d$ (it corresponds to a~$(D+1-d)$-bubble) and so its dimension is $d-1$.

\item In fact, since $\Delta$ is \emph{non-branching}, \emph{strongly connected} and \emph{pure}, it is a~\textit{$D$-dimensional simplicial pseudomanifold}~\cite{lost}.
\end{itemize}

\begin{Remark}[duality]\label{rem:complex}
The vertices of the graph correspond to the $D$-simplices of the simplicial complex. The half-lines of a vertex represent the $(D-1)$-simplices bounding a $D$-simplex and have a color. Any lower-dimensional sub-simplex is colored by the colors of the $D-1$ simplices sharing it. In Fig.~\ref{fig:complex} we sketched the dual complex in $D=3$ dimensions. The vertices are dual to tetrahedra. A triangle (say $3$) is dual to a line (of color~$3$) and separates two tetrahedra. An edge (say common to the triangles $2$ and $3$) is dual to a face (2-bubble of colors $2$ and~$3$). A~vertex (say common to the triangles $0$, $2$ and~$3$) is dual to a $3$-bubble (of colors $0$, $2$ and $3$).
\begin{figure}[t]\centering
\includegraphics[scale = 1.4]{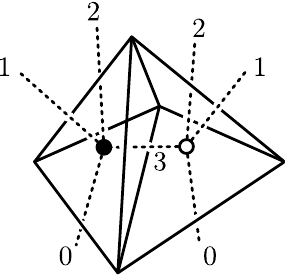}
\caption{The dual complex in $D=3$.}\label{fig:complex}
\end{figure}
\end{Remark}

\begin{Remark}[boundary graph]\label{rem:boundary-dual}
Let us return for a moment to examine the structure of a~boun\-dary graph ${\mathcal G}_{\partial}$. Such a~graph may have multiple connected components but each component has a cellular complex structure. For $d\ge 1$ the {boundary $d$-bubbles}~$({\mathcal B}_{\partial})^{i_1\dots i_{d+1}}_{(\rho)}$ are the maximally connected components of ${\mathcal G}_\partial $ formed by boundary vertices $v^{i_a}$ and boundary edges~$l^{i_bi_c}$, where $i_a,i_b, i_c \in \{i_1,\dots, i_{d+1}\}$. Following an analogous construction to that outlined above, but taking into account that the boundary $d$-bubbles have $d+1$ colors, one can show that each connected component of ${\mathcal G}_\partial$ is dual to a~simplicial complex (and is a~pseudomanifold). In fact, the simplicial complex dual to~${\mathcal G}_\partial$ is the boundary of the simplicial complex dual to~${\mathcal G}$.

Consequently, this motivates the study of these structures homological properties.
\end{Remark}

\subsection{Colored homology}

The topological spaces def\/ined by $(D+1)$-colored graphs are amenable to both a homological and homotopical analysis. The colored homology is def\/ined for the graph complex \cite{color} as follows:
\begin{Definition}[chain group] The {$d$-th chain group} is the group f\/initely generated by the $d$-bubbles
\begin{gather*}
 \alpha_d=\sum_{{\cal B}^{i_1 \dots i_d }_{(\rho)} \subset {\mathcal G} } c^{i_1 \dots i_d }_{(\rho)}
 {\cal B}^{i_1 \dots i_d }_{(\rho)} , \qquad c^{i_1 \dots i_d }_{(\rho)} \in \mathbb{Z} .
 \end{gather*}
\end{Definition}

The chain groups def\/ine homology groups via a boundary operator:
\begin{Definition}[boundary operator]\label{def:bound}
The {$d$-th boundary operator $\partial_d$} acting on a $d$-bubble ${\cal B}^{i_1 \dots i_d }_{(\rho)}$ is:
\begin{itemize}\itemsep=0pt
\item For $d\ge 2$,
 \begin{gather*} 
\partial_d( {\cal B}^{i_1 \dots i_d }_{(\rho)} ) = \sum_{q=1}^d (-1)^{q+1}
\sum_{ {\mathcal B}^{i_1 \dots \widehat{i_q} \dots i_d }_{(\kappa)} \subset {\cal B}^{i_1 \dots i_d }_{(\rho)} }
 {\mathcal B}^{i_1 \dots \widehat{i_q} \dots i_d }_{(\kappa)} ,
 \end{gather*}
which associates to a $d$-bubble the alternating sum of all $(d-1)$-bubbles formed by subsets of its vertices.

\item For $d=1$, since the edges ${\mathcal B}^i_{(\rho)}$ connect a positive vertex~$v$ to a negative one~$\bar v$:
 \begin{gather*}
 \partial_d {\mathcal B}^i_{(\rho)}= v - \bar v .
\end{gather*}
 \item For $d=0$, $\partial_0 v = \partial_0 \bar v =0$.
\end{itemize}
\end{Definition}

The colored boundary operators extend linearly over chains and thereafter def\/ine a homology as $\partial_{d-1}\circ \partial_{d}=0$ \cite{color}. Thus we def\/ine the {\it $d$-th colored homology group} to be $H_d \equiv \operatorname{ker}(\partial_d)/\operatorname{Im}(\partial_{d+1})$.

These graphs also facilitate a f\/inite presentation of their \textit{fundamental group} by associating a~generator to all edges of~${\mathcal G}$ (apart from those edges lying on a maximal tree) and a relation to all faces of~${\mathcal G}$.

\subsection{Equivalence: combinatorial and topological}

Colored graphs support a class of moves, termed $k$-dipole moves \cite{Gur3,Gur4,GurRiv}, that have a well controlled ef\/fect on their bubble structure. A~priori, these are combinatorial in nature and allow one to set up the notion of combinatorial equivalence.
However, a subset are homeomorphisms of the graph complex, thus setting the stage for topological equivalence.

\begin{Definition}[$k$-dipole] A {$k$-dipole} $d_k$ is a subset of ${\mathcal G}$ comprising of two vertices $v$, $\bar v$ such that:
\begin{itemize}\itemsep=0pt
\item $v$ and $\bar v$ share $k$ edges colored by $i_1, \dots, i_k\in {\mathbb Z}_{D+1}$;
\item $v$ and $\bar v$ lie in distinct $(D+1-k)$-bubbles: $B^{\hat{i}_1\dots \hat{i}_{k}}_{(\alpha)}\neq B^{\hat{i}_1\dots
\hat{i}_{k}}_{(\beta)}$.
\end{itemize}
 \end{Definition}

We say that $d_k$ \emph{separates} the bubbles $B^{\hat{i}_1\dots \hat{i}_{k}}_{(\alpha)}$ and $B^{\hat{i}_1\dots \hat{i}_{k}}_{(\beta)}$. Yet more important is how one manipulates the graph structure with respect to these subsets.

\begin{figure}[t]\centering
\includegraphics[width= 8cm]{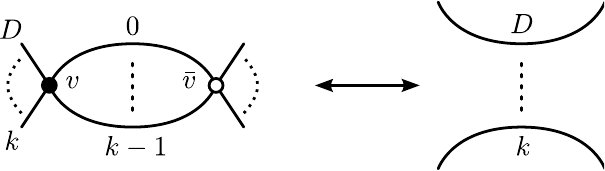}
\caption{The process of $k$-dipole contraction/creation in the case $i_p = p-1$ for all $p \in \{1,\dots, k\}$.}\label{fig:1canc}
\end{figure}

\begin{Definition}[$k$-dipole moves]
 The process of {$k$-dipole contraction}:
\begin{itemize}\itemsep=0pt
\item deletes the vertices $v_1$ and $v_2$;
\item deletes the edges $i_1,\dots, i_k$;
\item connects the remaining edges {\it respecting coloring}, see Fig.~\ref{fig:1canc}.
\end{itemize}
The process of $k$-dipole creation is precisely the inverse. We denote by ${\mathcal G}/d_k$ the graph obtained from ${\mathcal G}$ by contracting~$d_k$. Note that the separation property makes the identif\/ication and creation of $k$-dipoles somewhat subtle.
\end{Definition}

\begin{Definition}[combinatorial equivalence]\label{def:combeq}
Two graphs as said to be {combinatorially equi\-va\-lent}, denoted $\sim^{(c)}$, if they are related by a sequence of $k$-dipole contractions and creations.
 \end{Definition}

Together with combinatorial equivalence, comes the idea of combinatorial equivalence class. With the variety of moves available, these classes may have more or less restrictive membership requirements. One particularly useful set of classes respects 1-dipole moves. These are known as the \emph{combinatorial core equivalence classes}.
 \begin{Definition}[combinatorial core graph]\label{def:corecomb} A {combinatorial core graph at order $p$}, deno\-ted~${\mathcal G}^{(c)}_p$, is a~$(D+1)$-colored graph with $2p$ vertices, such that for all colors $i$, it has a unique $D$-bubble~${\mathcal R}^{\widehat{i}}_{(1)}$.
\end{Definition}

\textbf{Combinatorial bubble routing algorithm.} What is more, every colored graph is combinatorially equivalent to (at least) one combinatorial core graph via the contraction of a maximal set of 1-dipoles:
\begin{itemize}\itemsep=0pt
\item {\it Designate a root.} For a given color $i$, we pick one of the $D$-bubbles ${\mathcal B}^{\widehat{i}}_{(\rho)}$ as a root ${\mathcal R}^{\widehat{i}}_{(1)}$ bubble. The total number of roots of a graph is ${\mathcal R}^{[D]} = D+1$.

\item {\it Identify $\widehat i$-connectivity graph.} We associate to the bubbles $\widehat{i}$ of ${\mathcal G}$ an $\widehat{i}$-\emph{connectivity graph}. Its vertices represent the various bubbles ${\mathcal B}^{\widehat{i}}_{(\rho)} $. Its lines are the lines of color $i$ in ${\mathcal G}$. They either start and end on the same bubble $ {\mathcal B}^{\widehat{i}}_{(\alpha)} $, known as {\it tadpole lines} in the connectivity graph, or they do not. A particularly simple way to picture the $\widehat{i}$-connectivity graph is to draw ${\mathcal G}$ with the lines $j\neq i $ much shorter than the lines $i$.

\item {\it Choose a tree.} We chose a tree ${\mathcal T}^{i}$ in the $\widehat{i}$-connectivity graph, such that its root is ${\mathcal R}^{\widehat{i}}_{(1)}$. We refer to the rest of the lines of color $i$ as {\it loop lines}.

\item {\it Contract.} All the ${\mathcal B}^{[\widehat{i}]} - 1 $ lines of ${\mathcal T}^i$ are 1-dipoles and we contract them. We end up with a connectivity graph with only one vertex corresponding to the root bubble~${\mathcal R}^{\widehat{i}}_{(1)}$. The remaining lines of color $i$ cannot be contracted further (they are tadpole lines in the connectivity graph). The number of the $D$-bubbles of the other colors is unchanged under these contractions.

\item {\it Repeat.} We iterate the previous three points for all colors starting with $D$. The routing tree ${\mathcal T}^{j}$ is chosen in the graph obtained {\it after} contracting ${\mathcal T}^{j+1}, \dots, {\mathcal T}^D $. The number of bubbles of colors $q>j$ are constant under contractions of $1$-dipoles of color~$j$, hence the latter {\it cannot} create new 1-dipoles of color~$q$. Reducing a full set of 1-dipoles indexed by $D+1$ routing trees ${\mathcal T}^0, \dots, {\mathcal T}^{D} $ we obtain a graph in which all bubbles are roots. This is precisely a \emph{combinatorial core graph}.
\end{itemize}
This gives rise to:
\begin{Definition}[combinatorial core equivalence class] A combinatorial core equivalence class is a set of graphs related by sequences of 1-dipole moves.
\end{Definition}
\begin{Remark}[non-uniqueness of rooting] The rooting algorithm allows us to pick a representative for each core equivalence class. However, the combinatorial core graph one obtains by the above routing procedure is \emph{not} independent of the routing trees. The same graph leads to several equivalent core graphs, all at the same order~$p$.
\end{Remark}

We now brief\/ly mention the topological analogue. We shall utilize a fundamental result from combinatorial topology \cite{FG,Lins}:

\begin{Theorem}[topological $k$-dipole] Two pseudomanifolds dual to ${\mathcal G}$ and ${\mathcal G}/d_k$ are homeomorphic if one of the $(D+1-k)$-bubbles ${\mathcal B}^{\widehat{i}_0\dots \widehat{i}_{k-1}}_{(\alpha)} $ or ${\mathcal B}^{\widehat{i}_0\dots \widehat{i}_{k-1} }_{ (\beta) }$ separated by the dipole is dual to a~sphere~$S^{D-k}$.
\end{Theorem}
\noindent This allows us to propose another equivalence relation on the set of colored graphs:
\begin{Definition}[topological equivalence]
Two graphs are said to be {topologically equivalent}, denoted $\sim^{(t)}$, if they are
related by a sequence of topological dipole contraction and creation moves.
\end{Definition}
With appropriate modif\/ications, one can def\/ine topological core graph, topological core equivalence class and topological bubble rooting algorithm. We relinquish the details to other sources~\cite{GurRiv, razvan-review}.

Moreover, $(D+1)$-colored graphs representing manifolds are related, through the combinatorial/topological bubble routing algorithms to core graphs that are called crystallizations in the graph-encoded manifold literature~\cite{ferri-review}. The combinatorial and topological bubble routing algorithms coincide for manifolds, since all $D$-bubbles are homeomorphic to spheres and the routing algorithms only involve dipole contraction\footnote{Combinatorial 1-dipole insertion may produce two non-spherical $D$-bubbles from an initial spherical one.}.

\subsection{Jackets and degree}

As we have seen, the bubble structure captures the rich topology supported by these colored graphs. As a result, it tends to be rather intricate and subtle. We should also like to have a~somewhat blunter tool that captures only some of the information encoded by the colors.

This tool is provided by the \emph{jackets}. Their main advantage is that they are just ribbon graphs, like those generated by matrix models. As such they are Riemann surfaces embedded in the cellular complex and thus the subset of topological information that they capture is nicely encapsulated by their genera.

\begin{Definition}[jacket] A colored {jacket} ${\mathcal J}$ is a 2-subcomplex of ${\mathcal G}$, labeled by a $(D+1)$-cycle~$\tau$, such that:
 \begin{itemize}\itemsep=0pt
 \item ${\mathcal J}$ and ${\mathcal G}$ have identical vertex sets, ${\mathcal V}_{\mathcal J} = {\mathcal V}_{\mathcal G}$;
 \item ${\mathcal J}$ and ${\mathcal G}$ have identical edge sets, ${\mathcal E}_{\mathcal J} = {\mathcal E}_{\mathcal G}$;
 \item the face set of ${\mathcal J}$ is a subset of the face set of ${\mathcal G}$: ${\mathcal F}_{{\mathcal J}} = \{f\in {\mathcal F}_{{\mathcal G}}\,|\, f = (\tau^q(0),\tau^{q+1}(0)),\, q\in {\mathbb Z}_{D+1} \}$.
 \end{itemize}
 \end{Definition}
 \begin{Remark}[connectivity and abundance] It is evident that ${\mathcal J}$ and ${\mathcal G}$ have the same connectivity. In actual fact, a given jacket is independent of the overall orientation of the cycle, meaning that the number of jackets is in one-to-two correspondence with $(D+1)$-cycles. Therefore, the number of independent jackets is $D!/2$ and the number of jackets containing a given face is $(D-1)!$.
\end{Remark}

\begin{Remark}[jackets as ribbon graphs] The jacket has the structure of a \emph{ribbon graph}. Note that each edge of ${\mathcal J}$ lies on the boundary of two of its faces. Thus, it corresponds to a ribbon line in the ribbon graph. As we said, the ribbon lines separate two faces, $(\tau^{-1}(i),i)$
and $(i,\tau(i))$ and inherit the color $i$ of the line in ${\mathcal J}$. Ribbon graphs are well-known to correspond to Riemann surfaces, and so the same holds for jackets. Given this, we can def\/ine the \emph{Euler characteristic} of the jacket as: $\chi({\mathcal J}) = |{\mathcal F}_{\mathcal J}| - |{\mathcal E}_{\mathcal J}| + |{\mathcal V}_{\mathcal J}| = 2 - 2g_{\mathcal J}$, where $g_{{\mathcal J}}$ is the \emph{genus} of the jacket\footnote{A moment's ref\/lection reveals that the jackets necessarily represent orientable surfaces.}.
\end{Remark}

\begin{Remark}[examples] In $D=2$, the (unique) jacket of a $(2+1)$-colored graph is the graph itself. An example of a graph and its jackets (and their associated cycles) is given in Fig.~\ref{fig:exemplujacket1}. For instance the leftmost jacket corresponding to the cycle $\tau = (0123)$ contains only the faces~$01$, $12$, $23$ and $30$.
\begin{figure}[t]\centering
\includegraphics[width=12cm]{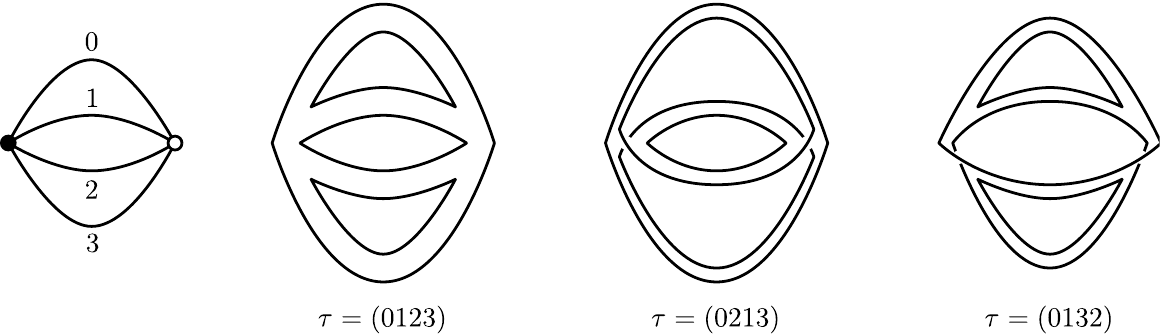}
\caption{A vacuum 4-colored graph and its jackets.}\label{fig:exemplujacket1}
\end{figure}
\end{Remark}

\begin{Remark}[jackets and $d$-bubbles] For a $(D+1)$-colored graph ${\mathcal G}$, its $D$-bubbles are $D$-colored graphs ${\mathcal B}_{(\rho)}^{\widehat i}$. Thus, they also possess jackets, which we denote by ${\mathcal J}_{(\rho)}^{\widehat i}$. It is rather elementary to construct the ${\mathcal J}^{\widehat{i}}_{(\rho)}$
 from the ${\mathcal J}$. Let us construct the ribbon graph ${\mathcal J}^{\widehat i}$ consisting of vertex, edge and face sets:
\begin{gather*}
{\mathcal V}_{{\mathcal J}^{\widehat i}} = {\mathcal V}_{{\mathcal J}}, \qquad {\mathcal E}_{{\mathcal J}^{\widehat i}} = {\mathcal E}_{{\mathcal J}}\setminus E^i, \qquad \textrm{and}\\
 {\mathcal F}_{{\mathcal J}^{\widehat{i}}}
 = \big({\mathcal F}_{\mathcal J} \setminus \big\{\big(\tau^{-1}(i), i\big), (i,\tau(i))\big\}\big) \cup \big\{ \big(\tau^{-1}(i), \tau(i)\big) \big\} ,
\end{gather*}
that is having all vertices of ${\mathcal G}$, all lines of ${\mathcal G}$ of colors dif\/ferent from $i$ and some faces. Given that the face set of ${\mathcal J}$ is specif\/ied by a $(D+1)$-cycle $\tau$, the f\/irst thing to notice is that the face set of ${\mathcal J}^{\hat i}$ is specif\/ied by a $D$-cycle obtained from~$\tau$ by deleting the color~$i$. The ribbon sub\-graph~${\mathcal J}^{\widehat{i}}$ is the union of several connected components, ${\mathcal J}^{\widehat{i}}_{(\rho)}$. Each~${\mathcal J}^{\widehat{i}}_{(\rho)}$ is a jacket of the $D$-bubble~${\mathcal B}^{ \widehat{i} }_{ (\rho) }$. Conversely, every jacket of ${\mathcal B}^{ \widehat{i} }_{ (\rho) }$ is obtained from exactly $D$ jackets of~${\mathcal G}$. To realize this, consider a jacket ${\mathcal J}^{\widehat{i}}_{(\rho)}$. It is specif\/ied by a~$D$-cycle (missing the color $i$). On can insert the color $i$ anywhere along the cycle and thus get $D$ independent $(D+1)$-cycles.

More generally, the $d$-bubbles are $d$-colored graphs and they also possess jackets which can be obtained from the jackets of~${\mathcal G}$.

Consider once again Fig.~\ref{fig:exemplujacket1}. Applying our procedure to the jacket $(0123)$ leads to the three jackets $(123)$, $(023)$ and $(012)$. Each of these jackets corresponds to a bubble of Fig.~\ref{fig:exemplubub} and is a~3-colored graph.
\end{Remark}

\begin{Definition}[degree] We def\/ine:
 \begin{itemize}\itemsep=0pt
 \item the (convergence) degree of a graph ${\mathcal G}$ is $\omega({\mathcal G})=\sum_{{\mathcal J}} g_{{\mathcal J}}$, where the sum runs over all the jackets ${\mathcal J}$ of ${\mathcal G}$,
 \item the degree of a $k$-dipole $d_k$ is the lesser of the degrees of the two $(D+1-k)$-bubbles that it separates.
 \end{itemize}
\end{Definition}

\begin{Remark}[properties of the degree] The degree displays a number of pertinent properties:
\begin{itemize}\itemsep=0pt
 \item The degree is a non-negative integer that is readily computable from the graph (as it is the sum of the genera of the embedded jackets).

 \item The degree of a graph $\omega({\mathcal G})$ and its $D$-bubbles $\omega\big({\mathcal B}^{\widehat{i}}_{(\rho)}\big)$ are not independent, that is, they respect
 \begin{gather*}
 \omega({\mathcal G}) = \frac{(D-1)!}{2} \big( p+D-{\mathcal B}^{[D]} \big) + \sum_{i;\rho} \omega\big({\mathcal B}^{\widehat{i}}_{(\rho)}\big) ,
 \end{gather*}
 where $2p$ is the number of vertices in ${\mathcal G}$ and ${\mathcal B}^{[D]}$ is the total number of $D$-bubbles of all colors.

\item The degree of a graph changes under $k$-dipole contraction ${\mathcal G} \rightarrow {\mathcal G}/d_k$. The degree of ${\mathcal G}$ and~${\mathcal G}/d_k$ are related by
 \begin{gather*}
 \omega({\mathcal G}) = \frac{(D-1)!} {2} \big( (D+1)k - k^2 - D \big) + \omega({\mathcal G}/d_k) .
 \end{gather*}
In particular, the degree is unchanged by any 1-dipole move. As a consequence, all graphs in a combinatorial core equivalence class have the same degree.

 \item Let ${\mathcal G}$ be a $D+1$ colored graph and ${\mathcal B}^{ \widehat{D} }_{(\rho)}$ its $D$-bubbles with colors $\widehat{D}$. Then
\begin{gather*}
 \omega({\mathcal G}) \ge D \sum_{\rho} \omega \big({\mathcal B}^{\widehat{D} }_{(\rho)} \big) .
\end{gather*}
\end{itemize}
\end{Remark}

\begin{figure}[t] \centering
 \includegraphics[scale = 1]{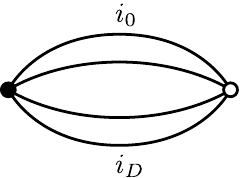}
 \qquad \quad
 \includegraphics[scale = 1]{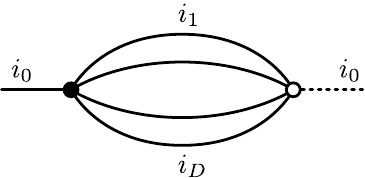}
 \caption{The supermelon (left) and an elementary melon (right).}\label{fig:melon}
\end{figure}

\begin{Remark}[ramif\/ications of vanishing degree: melonic graphs] \label{rem:melonic} One particular core equivalence class is distinguished by vanishing degree $\omega({\mathcal G}) = 0$. Of course, the fact that all graphs with vanishing degree lie in the same core equivalence class must be shown explicitly. This cannot be taken for granted.
\begin{itemize}\itemsep=0pt
\item A simple argument yields that $p+D-{\mathcal B}^{[D]} \geq 0$.\footnote{The quantity $p - {\mathcal B}^{[D]}$ is conserved by 1-dipole moves. A core graph has at least two vertices ($p \geq 1$) and ${\mathcal B}^{[D]} = D+1$. Thus, $p+D- {\mathcal B}^{[D]} \geq 0$ for core graphs. Since every graph roots to a~core graphs, the result follows.} Thus, one has that $\omega({\mathcal G}) = 0 \implies \omega({\mathcal B}^{\hat{i}}_{(\rho)}) = 0$ for all $i$ and all~$\rho$, as well as $p +D - {\mathcal B}^{[D]} = 0$.

\item Since a core graph has precisely ${\mathcal B}^{[D]} = D+1$, one with vanishing degree must have $p = 1$. There is a unique core graph with~2 vertices, known af\/fectionately as the \emph{supermelon}. It is illustrated in Fig.~\ref{fig:melon}.

\item In $D = 2$, the 1-dipole moves allow one to explore the full set of colored planar graphs. Our interested lies with higher values of $D$, however, so we shall not explore this further.

\item In $D \geq 3$, $\omega({\mathcal G}) = 0$ implies that ${\mathcal G}$ contains an \emph{elementary melon}, illustrated in Fig.~\ref{fig:melon}. Such a sub-graph can be removed trivially through 1-dipole reduction. Iterating this process on the resulting graph means that all graphs with vanishing degree root to the supermelon through the removal of a sequence of elementary melons. Hence this core equivalence class is referred to as the class of melonic graphs.

\item
 If $\omega({\mathcal G})=0$ then ${\mathcal G}$ is dual to a sphere $S^D$. The reciprocal holds in $D=2$. Thus, in all melonic graphs are spheres.
 \end{itemize}
\end{Remark}

\begin{Remark}[some further results in $D = 3$]
Let us focus brief\/ly on the case when $D = 3$, where we can readily utilize this machinery to uncover some further results about the topology of the 4-colored graphs.
\begin{itemize}\itemsep=0pt
 \item The jacket structure can tell us more about the topology of the full graph, namely: \emph{if ${\mathcal G}$ possesses a spherical jacket then~${\mathcal G}$ is spherical}.

 \item The jackets have a special signif\/icance as splitting surfaces~\cite{heegaard-tensor}: \emph{if ${\mathcal G}$ is a manifold, then its jackets~${\mathcal J}$ are Heegaard surfaces}, where a \emph{Heegaard splitting} of a compact connected oriented 3-manifold~${\mathcal M}$ is an ordered triple $(\Sigma, {\mathcal H}_1, {\mathcal H}_2)_{\mathcal M}$ consisting of a~compact connected oriented surface $\Sigma$ and two handlebodies ${\mathcal H}_1$, ${\mathcal H}_2$ such that $\partial{\mathcal H}_1 = \partial{\mathcal H}_2 = \Sigma$. $\Sigma$ is known as the \emph{Heegaard surface} of the splitting.
\end{itemize}
\end{Remark}

\section{Algebraic properties}\label{sec:alg}

Colored graphs support a very nuanced algebraic structure \cite{alg1,alg2}. Its importance for the understanding of tensor models cannot be overstressed. In principle, it provides a bridge between the apparent discrete world inhabited by the tensor representation and any potential continuum representation. In reality, work on this facet of the theory has gone little further than its def\/inition and the demonstration of its self-consistency.

To begin, we need some basic graph-theoretic concepts:
\begin{Definition}[marked graph]
 A marked graph $({\mathcal B}, \bar{v})$ is a $D$-colored graph ${\mathcal B}$ paired together with one of its negative vertices $\bar{v}$.
\end{Definition}
\begin{Definition}[graph contraction]
 Consider two $D$-colored graphs ${\mathcal B}_1$ and ${\mathcal B}_2$ and a positive-negative vertex pair $(v_1,\bar{v}_2)\in{\mathcal B}_1\times{\mathcal B}_2$ (or $(\bar{v}_1,v_2)$). Contraction of ${\mathcal B}_1$ and ${\mathcal B}_2$ at $(v_1,\bar{v}_2)$ is obtained in a two-step process:
 \begin{itemize}\itemsep=0pt
 \item delete $v_1$ and $\bar{v}_2$ along with the half-edges emanating from them;
 \item reconnect the surviving half-edges in the (unique) manner that preserves the color structure.
 \end{itemize}
 One denotes the resulting $D$-colored graph by ${\mathcal B}_1\star_{(v_1,\bar{v}_2)}{\mathcal B}_2$. Note that this def\/inition coincides with that of $0$-dipole contraction for $D$-colored graphs.
\end{Definition}
\begin{Remark}[property] Note also that graph contraction has one vital property
 \begin{gather} \label{eq:contraction-property}
 ({\mathcal B}_1\star_{(v,\bar{v}_2)}{\mathcal B}_2)\star_{(w,\bar{v}_3)}{\mathcal B}_3 = {\mathcal B}_1\star_{(v,\bar{v}_2)}({\mathcal B}_2\star_{(w,\bar{v}_3)}{\mathcal B}_3).
 \end{gather}
\end{Remark}

 \textbf{Constructing the bubble algebra:}

\begin{itemize}\itemsep=0pt
 \item Consider the set $S = \{{\mathcal L}_{({\mathcal B},\bar{v})}\}$ of elements indexed by marked $D$-colored graphs.

\item Construct an (inf\/inite-dimensional) vector space over the reals using the elements of $S$ as basis vectors. Denote this vector space by~$X$.

\item Endow $X$ with the following non-associative multiplication
 \begin{gather*} 
 \big[{\mathcal L}_{({\mathcal B}_1,\bar{v}_1)},{\mathcal L}_{({\mathcal B}_2,\bar{v}_2)}\big] \equiv
 \sum_{v\in{\mathcal B}_1}{\mathcal L}_{({\mathcal B}_1\star_{(v,\bar{v}_2)}{\mathcal B}_2,\bar{v}_1)} -
 \sum_{v\in{\mathcal B}_2}{\mathcal L}_{({\mathcal B}_2\star_{(v,\bar{v}_1)}{\mathcal B}_1,\bar{v}_2)}.
 \end{gather*}
\item Impose \emph{bilinearity} in the f\/irst argument by f\/iat
 \begin{gather*}
 \big[a{\mathcal L}_{({\mathcal B}_1,\bar{v}_1)} + b{\mathcal L}_{({\mathcal B}_2,\bar{v}_2)}, {\mathcal L}_{({\mathcal B},\bar{v})}\big] \equiv a\big[{\mathcal L}_{({\mathcal B}_1,\bar{v}_1)}, {\mathcal L}_{({\mathcal B},\bar{v})}\big] + b\big[{\mathcal L}_{({\mathcal B}_2,\bar{v}_2)}, {\mathcal L}_{({\mathcal B},\bar{v})}\big]
 \end{gather*}
 along with a similar relation for the second argument.
\item This bracket is clearly \emph{anticommutative}, since
 \begin{gather*} 
 \big[{\mathcal L}_{({\mathcal B}_1,\bar{v}_1)},{\mathcal L}_{({\mathcal B}_2,\bar{v}_2)}\big] = -\big[{\mathcal L}_{({\mathcal B}_2,\bar{v}_2)},{\mathcal L}_{({\mathcal B}_1,\bar{v}_1)}\big],
 \end{gather*}
 and this extends to all of $X$ using bilinearity of the bracket.

 In fact, with the help of \eqref{eq:contraction-property}, the bracket satisf\/ies the \emph{Jacobi identity}
 \begin{gather*}
 \big[\big[{\mathcal L}_{({\mathcal B}_1,\bar{v}_1)},{\mathcal L}_{({\mathcal B}_2,\bar{v}_2)}\big],{\mathcal L}_{({\mathcal B}_3,\bar{v}_3)}\big]
 + \big[\big[{\mathcal L}_{({\mathcal B}_2,\bar{v}_2)},{\mathcal L}_{({\mathcal B}_3,\bar{v}_3)}\big],{\mathcal L}_{({\mathcal B}_1,\bar{v}_1)}\big] \\
\qquad {}+ \big[\big[{\mathcal L}_{({\mathcal B}_3,\bar{v}_3)},{\mathcal L}_{({\mathcal B}_1,\bar{v}_1)}\big],{\mathcal L}_{({\mathcal B}_2,\bar{v}_2)}\big] = 0.
\end{gather*}
\end{itemize}
In the end, we have the following neat result:
\begin{Theorem}[bubble algebra] The pair $(X,[\cdot,\cdot])$ form a~Lie algebra.
\end{Theorem}

The next step is to look for representations of this algebra. One such is given by the tensor model $D$-bubble observables of a $(D+1)$-colored tensor model, hence the rather nomenclature. This algebra underpins the symmetries of the tensor model theory.

Important for the representation theory is the identif\/ication of subalgebras, from which one can induce representations of the full algebra. In this direction, one has:
\begin{Theorem}[melonic subalgebra] The marked melonic $D$-colored graphs form a subalgebra of the $D$-colored bubble algebra.
\end{Theorem}

\section[Melonic graphs: probing deeper into their combinatorial and metric structures]{Melonic graphs: probing deeper into their combinatorial\\ and metric structures}\label{sec:mel}

In order to perform a more detailed analysis, let us restrict ourselves to a single core equivalence class: the melonic graphs that we introduced earlier. This is the sole core equivalence class contributing to the leading order of generic tensor model theories and so a thorough understanding of its combinatorial and metric properties is of utmost importance.

\looseness=-1 By def\/inition, melonic graphs are those that reduce to the \emph{supermelon} core graph. Thus, they have vanishing degree and dif\/fer from the supermelon by a sequence of 1-dipole moves. But this is a rather redundant prescription as many sequences lead to the same melonic graph. However, as laid out in Remark \ref{rem:melonic}, melonic graphs reduce to the supermelon through the iterative removal of \emph{elementary melons}. This provides a parsimonious prescription for melonic graphs.

{\bf Rooted melonic graphs.} To be precise, we shall construct \emph{rooted melonic graphs}. A~rooted melonic graph is a melonic graph with one edge singled out. Such rooting is common in combinatorial graph theory, as it simplif\/ies counting problems. Marking an edge of the supermelon graph essentially yields an elementary melon of some color~$i$ (the color of the distinguished edge). In the following we shall always root the graphs along an edge of color~$0$.

The elementary melons also act as the fundamental building blocks of generic melonic graphs. To add a bit more nomenclature, an elementary melon consists of two vertices connected by~$D$ edges. Both vertices have one external edge. Obviously, both external edges possess the same color, say~$i$. An elementary melon has two features: \textit{i})~an external edge of color~$i$ incident to the white vertex, which is known as
the \emph{root edge}; \textit{ii})~$D+1$ edges incident at the black vertex, which are known as \emph{active edges}, having distinct colors from $\{0,1,\dots, D\}$.

{\bf Rooted melonic graph construction algorithm.} One can construct the class rooted melonic graphs iteratively.
\begin{description}\itemsep=0pt
\item[$p=1$:] There is a unique rooted melonic graph with two vertices. It is illustrated in the bottom left of Fig.~\ref{fig:elMelonInsertion}
and is the \emph{elementary melon of color~$0$}.
\begin{figure}[t]\centering
\includegraphics[scale=1]{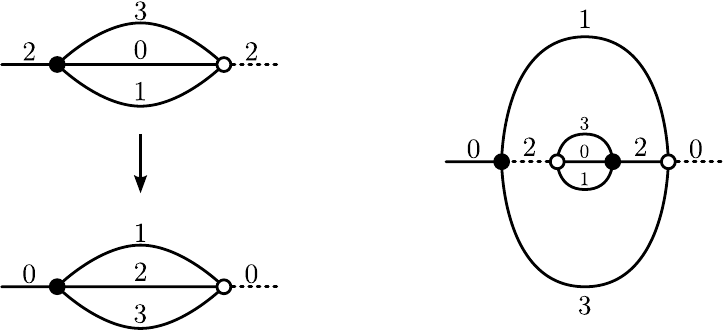}
\caption{An elementary melon of color $2$ inserted along the active edge of color $2$ (for $D=3$). The active
edges are drawn using full lines.}\label{fig:elMelonInsertion}
\end{figure}

\item[$p=2$:] There are $D+1$ melonic graphs with four vertices. One obtains them from the graph at $p=1$ by replacing an active edge of a given color by an elementary melon of the same color (as shown in Fig.~\ref{fig:elMelonInsertion}).

\item[$p=k$:] One obtains these graphs from those at $p=k-1$ by replacing some active edge by an elementary melon of the appropriate color.
\end{description}

The need for a precise prescription cannot be overstressed, since all the ensemble properties of these graphs, stem from being able to count them precisely.

As mentioned earlier, the abstract structure of rooted melonic graphs coincides with that of several other objects, which we shall describe presently.

{\bf Rooted melonic graphs as colored rooted $\boldsymbol{(D+1)}$-ary trees.} There is a simple bijection between the set of rooted melonic graphs and \emph{colored rooted $(D+1)$-ary trees}. The fundamental building blocks of any colored rooted $(D+1)$-ary tree are the \emph{elementary vertices}. An~elementary vertex of color $i$ is $(D+2)$-valent with two distinguished features: \textit{i})~a \emph{root edge} of color~$i$; \textit{ii})~$D+1$ \emph{active leaves} each with a distinct color from $\{0,\dots,D\}$. These correspond to the root edge and the active edges of the elementary melon, respectively. Since this class of trees is also constructed in an iterative manner, the map is self-evident:
\begin{description}\itemsep=0pt
\item[$p=1$:] There is a unique colored rooted $(D+1)$-ary tree with a single elementary vertex. This is the elementary vertex of color 0.
It is illustrated in the bottom left of Fig.~\ref{fig:elTree}.
\begin{figure}[t]\centering
\includegraphics[scale=1]{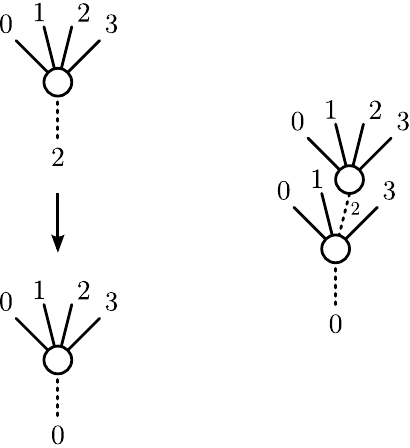}
\caption{The elementary vertex of color 2 replacing a leaf of color~2 (for $D=3$).}\label{fig:elTree}
\end{figure}

\item[$p=2$:] There are $D+1$ such trees with two elementary vertices. One obtains them from the tree at $p=1$ by replacing a leaf of a given color with an \emph{elementary vertex} of the same color (as shown in Fig.~\ref{fig:elTree}).

\item[$p=k$:] One obtains these trees from those at $p=k-1$ by replacing a leaf with an elementary vertex of the same color.
\end{description}

{\bf Rooted melonic graphs as colored simplicial $\boldsymbol{D}$-balls.} The description of how any given $(D+1)$-colored graph is dual, in a~precise topological sense, to a unique $D$-dimensional abstract simplicial pseudomanifold was provided earlier.

Consider cutting a closed $(D+1)$-colored graph along one edge. This results in an open graph whose dual is a simplicial complex with boundary. This boundary is a $(D-1)$-sphere constructed from two $(D-1)$-simplices.

Melonic graphs are dual to simplicial $D$-spheres. Rooted melonic graphs, which are melonic graphs with one edge cut, are dual to simplicial $D$-balls with the boundary just mentioned. For the want of a better name, we shall call them \emph{melonic $D$-balls}. One can def\/ine them iteratively. The fundamental building blocks are the \emph{elementary melonic $D$-balls}. These consist of two $D$-simplices sharing $D$ of their $(D-1)$-simplices. There are two more $(D-1)$-simplices forming the boundary $(D-1)$-sphere. In the manner outlined above, the two $D$-simplices are dual to the two vertices of the elementary melon, while the $(D-1)$-simplices are dual to the edges. Thus, the $(D-1)$-simplices inherit a single color. An elementary melonic $D$-ball of color $i$ has two distinguished features: \textit{i})~an external \emph{root $(D-1)$-simplex} of color $i$; \textit{ii})~ $D+1$ \emph{active $(D-1)$-simplices} (one of which is on the boundary), each with a~distinct color. The iterative def\/inition proceeds as follows:

\begin{description}\itemsep=0pt
\item[$p=1$:] There is a unique $D$-ball comprised of two $D$-simplices. It is the elementary melonic $D$-ball of color~0. It is illustrated in the bottom left of Fig.~\ref{fig:elballinsertion}.
\begin{figure}[t]\centering
\includegraphics[scale=1]{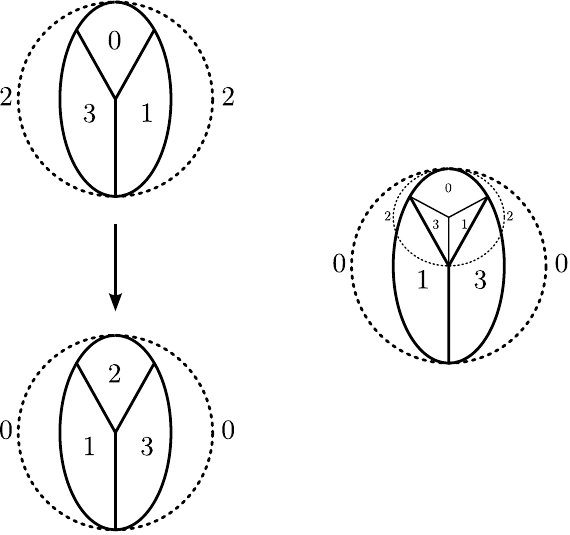}
\caption{The melonic $D$-ball at $p=2$ (for $D=3$), obtained by adding an elementary $D$-ball of color~2.}\label{fig:elballinsertion}
\end{figure}

\item[$p=2$:] There are $D+1$ melonic $D$-balls with four $D$-simplices. One obtains them from the melonic $D$-ball at $p=1$ by \textit{adding} an elementary melonic $D$-ball of a given color (shown in Fig.~\ref{fig:elballinsertion}). More precisely, one splits the melonic $D$-ball arising at $p=1$ along an active $(D-1)$-simplex (or one selects the active boundary $(D-1)$-simplex). One then glues an elementary melonic $D$-ball of the appropriate color along the split (or simply on the boundary $(D-1)$-simplex).

 \item[$p=k$:] One obtains them from those at $p=k-1$ by adding an elementary melonic $D$-ball at some active $(D-1)$-simplex.
\end{description}

\begin{Remark}[vertex correspondance] Let us denote a colored rooted $(D+1)$-ary tree with $p$ elementary vertices by $\mathcal{T}_p$ and its set of elementary vertices by $t_p$.

Noting that an elementary $D$-ball has precisely one \emph{internal} vertex before it is inserted, one notes that a melonic $D$-ball constructed form~$p$ elementary $D$-balls has precisely $p$ internal vertices. We shall denote such a melonic $D$-ball by $\mathcal{M}_p$ and its internal vertex set by~$m_p$.

Moreover, given such melonic $D$-ball $\mathcal{M}_p$, along with its associated $(D+1)$-ary tree~$\mathcal{T}_p$, one notes that their respective vertex sets are in bijective correspondence, since an elementary tree vertex corresponds to an elementary $D$-ball, which in turn has one internal vertex (before insertion).
\end{Remark}

{\bf Branches, words and vertex ordering.} Consider again an associated pair $\mathcal{M}_p$ and $\mathcal{T}_p$. As a tree, $\mathcal{T}_p$ has \emph{branches}, joining the \emph{root tree vertex} to each of the elementary vertices in~$t_p$. (Note that the root tree vertex is not in~$t_p$. Rather it is the other endpoint of the root edge in the f\/irst elementary tree vertex of color~$0$.)

\begin{figure}[t] \centering
\includegraphics[scale=1.2]{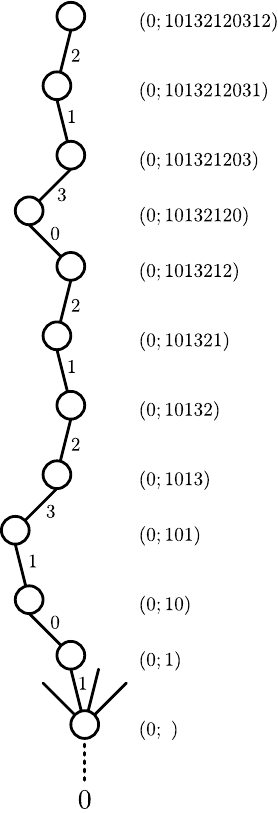}
\caption{The words associated to some vertices of a colored rooted $(D+1)$-ary tree.}\label{fig:words}
\end{figure}

\looseness=1 To each elementary tree vertex, one may associate a~\emph{word} drawn from the alphabet $\Sigma_{D+1} = \{0,1,\dots,D\}$. This word is constructed by listing left-to-right the colors encountered as one traverses the branch from the root vertex to the elementary vertex in question. As examples, the word associated to the root vertex is~$(0;\ )$, indicating the root edge of color $0$, while the word associated to the generic example given in Fig.~\ref{fig:words} is $(0;10132120312)$. In turn, these words have a natural (lexicographical) ordering, which may be used to convert~$t_p$ into an ordered set.

Due to the bijection relating $m_p$ to $t_p$, the internal vertices of the melonic $D$-ball are ordered. Moreover, each element of $m_p$ also inherits a word from its associated tree vertex. That word is not meaningless in this context. On the contrary, it captures the sequence of elementary $D$-balls, each one nested inside the preceding, on route to inserting that internal vertex.

\subsection[Generating rooted melonic $D$-balls and the continuum limit]{Generating rooted melonic $\boldsymbol{D}$-balls and the continuum limit}

Ultimately, f\/inite melonic $D$-balls hold only a limited amount of interest for us. We wish to know about the conf\/igurations in the large-$p$ limit, where one might conceivably formulate a~meaningful continuum limit~\cite{gur-leading}. Na\"{\i}vely, one could imagine that this limit is arrived at by considering conf\/igurations with increasingly large numbers of tetrahedra of diminishing individual size.

Our aim here is to make this initial idea more precise. In that respect, let us partition the set of melonic $D$-balls according to the number of internal vertices. At a given value of $p$, we want a~method to examine the properties, \emph{on average}, of that subset of melonic $D$-balls. A convenient method is to utilize a generating function approach
\begin{gather*}
 G(z) = \sum_{\mathcal{M}} z^{p_\mathcal{M}} = \sum_{p = 0}^\infty C^{(D+1)}_p z^p,\qquad \textrm{where}\quad C^{(D+1)}_p = \sum_{\mathcal{M}\colon p_{\mathcal{M}} = p} 1, \qquad z\in\mathbb{C}.
\end{gather*}
The f\/irst sum is over all melonic $D$-balls, weighted according to the number of internal vertices. Thus, the subset of melonic $D$-balls with $p$ internal vertices are all weighted equally allowing the reduction to the second sum, where $C^{(D+1)}_p$ is the cardinality of that subset.

One can expect that $G(z)$ has a f\/inite radius of convergence $z_c$. As a result, one can examine properties of the large-$p$ limit by tuning the \emph{coupling constant} $z$ to its \emph{critical value} $z_c$, as it is the large-$p$ coef\/f\/icients that determine the behaviour of~$G(z)$ in that region.

\begin{figure}[t]\centering
 \includegraphics[scale = 1.5]{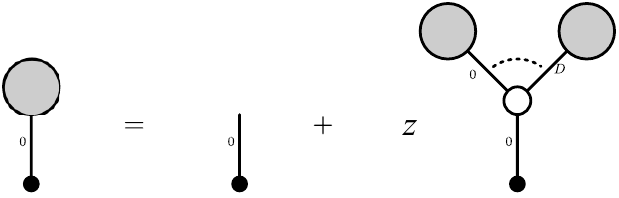}
 \caption{Consistency equation for the generating function~$G(z)$.}\label{fig:eq-trees}
\end{figure}

Thus, our f\/irst task is to gain control over the coef\/f\/icients $C^{(D+1)}_p$. In this regard, it is worth remembering that melonic $D$-balls are in a precise correspondence with colored rooted $(D+1)$-ary trees and experts have been counting trees since time immemoriam. A convenient way to count trees is to develop a consistency relation for the generating function $G(z)$. This is,
\begin{gather*}
 G(z) = 1 + z [G(z)]^{D+1},
\end{gather*}
which is illustrated in Fig.~\ref{fig:eq-trees}. It sums up the fact that the trees comprising $G(z)$ consist of the tree with just one root edge (and no elementary vertex) along with those with at least one elementary vertex (weighted by~$z$), along any whose~$D+1$ active edges one can insert an arbitrary colored rooted $(D+1)$-ary tree.

Expanding the terms and equating the coef\/f\/icients for the various powers of~$z$, one f\/inds a~recurrence relation for the $C^{(D+1)}_p$. This recurrence relation is satisf\/ied by the $(D+1)$-Catalan numbers
\begin{gather*}
C^{(D+1)}_p = \frac{1}{(D+1)p + 1}{(D+1)p + 1 \choose p}.
\end{gather*}
Applying Stirling's formula to the series coef\/f\/icients, one determines the large order behaviour
\begin{gather*}
 G(z) \sim \beta (z_c)^{-p}p^{-3/2} ,\qquad \textrm{where} \quad \beta = \dfrac{e}{\sqrt{2\pi}}\sqrt{\frac{D+1}{D^3}},\qquad z_c = \dfrac{D^D}{(D+1)^{D+1}}.
\end{gather*}
As expected, the above series has a f\/inite radius of convergence and the behaviour of the series in the vicinity of $z_c$ is
\begin{gather*}
 G(z) \sim \left(1 - \frac{z}{z_c}\right)^{1-\gamma} ,\qquad\textrm{where}\quad \gamma = \frac{1}{2}.
\end{gather*}
The exponent $\gamma$ is known as the \emph{susceptibility} and is the f\/irst example of a quantity that determines the properties of the continuum limit.

To this point, however, we have only analysed the large-$p$ limit. To obtain a continuum limit, we need to endow each melonic $D$-ball with a~metric. We shall develop that in more detail later. But for now, all we need is that under this metric, each melonic $D$-ball becomes an equilateral triangulation. This means that each tetrahedron has equal volume, known as the \emph{microscopic volume}~$\nu$. Meanwhile, for a melonic $D$-ball with $p$ internal vertices, and thus $2p$ tetrahedra, the \emph{macroscopic volume} is $V_\mathcal{M} = 2p\nu$. Thus, we may obtain a~f\/inite macroscopic volume if $\nu\rightarrow 0$ as $p\rightarrow \infty$ in some balanced fashion.

To do this meaningfully, remember that we are interested in the properties of melonic $D$-balls on average. Thus, it is worth analyzing
\begin{gather*}
 V(z) \equiv \sum_{\mathcal{M}} (2p_\mathcal{M}\nu)z^{p_{\mathcal{M}}} = \sum_{p = 0}^{\infty} (2p\nu) C^{(D+1)}_{p} z^{p} = 2\nu z\frac{\partial G}{\partial z} (z) \underset{z\sim z_c}{\sim} \nu z \left(1 - \frac{z}{z_c}\right)^{-1/2}.
\end{gather*}
If we consider taking the microscopic volume $\nu \rightarrow 0$ and $z \rightarrow z_c$ in a balanced fashion such that the macroscopic volume $V(z_c-)$ remains f\/inite. Then, we can realistically call such a~limit a~continuum limit.

\begin{Remark}[melonic $D$-balls and branched polymer spacetimes] In total, we shall extract three exponents associated with the continuum limit of melonic $D$-balls: the susceptibility \mbox{$\gamma = 1/2$}, the Hausdorf\/f dimension $d_H = 2$ and the spectral dimension $d_S = 4/3$ (although strictly speaking, the last refers to the dual rooted melonic graphs). These coincide with the exponents calculated for so-called \emph{branched polymer spacetimes} arising in the dynamical triangulations approach to quantum gravity~\cite{dt}. These have undesirable physical properties, among them the non-physical Hausdorf\/f and spectral dimensions. As a result, one concludes that melonic $D$-balls possess an uninteresting continuum limit from a physical spacetime viewpoint. However, these exponents are one of the major ways for testing the physical viability of any particular phase of a tensor model theory.
\end{Remark}

\subsection{Hausdorf\/f dimension}

The next exponent pertaining to the continuum limit of melonic $D$-balls is the Hausdorf\/f dimension:

\begin{Definition}[Hausdorf\/f dimension]
For a metric space $X$, the \emph{Hausdorff dimension} $d_H$ captures how the volume of a ball scales with respect to its geodesic radius.
More formally, it is def\/ined as
\begin{gather*}
d_H = \inf\{d\geq 0 \colon \mathcal{H}_d(X) = 0 \},
\end{gather*}
where $\mathcal{H}_d(X)$ is the $d$-dimensional Hausdorf\/f measure on $X$, that is
\begin{gather*}
\mathcal{H}_d(X) = \inf\bigg\{\delta = \sum_i r_i^{d}\colon \textrm{the indexed collection of balls of radius $r_i$ covers $X$}\bigg\}.
\end{gather*}
\end{Definition}

Obviously, for the simple example of f\/lat $D$-dimensional Euclidean space: $V_D \sim r^{D}$. Thus
\begin{gather*}
d_H = D.
\end{gather*}
It is clear, however, that some work must be done to extract this dimension for the class of melonic $D$-balls. One may follow the arguments in~\cite{Albenque, melondimension} for the full picture, which would be too laborious a task to reproduce in its entirety here. In fact, a detailed statement of the result already takes some considerable ink. So, with no further ado, the main result is:
\begin{Theorem}\label{th:GH} Under the uniform distribution, the family of melonic $D$-balls converges in the Gromov--Hausdorff topology on compact metric spaces to the continuum random tree
 \begin{gather*} 
 \left(m_p, \frac{d_{m_p}}{\Lambda_{\Delta} \sqrt{\frac{(D+1)p}{D}}} \right) \underset{p\to \infty}{\longrightarrow} ( {\cal T}_{2e}, d_{2e}) .
 \end{gather*}
\end{Theorem}
In turn, the Hausdorf\/f dimension of the limiting continuum metric space, the continuum random tree, is well known to be
\begin{gather*}
d_H = 2 ,
\end{gather*}
although this may be read of\/f directly from the above theorem, as we shall see presently.

{\bf Melonic $\boldsymbol{D}$-balls as metric spaces.} Consider a melonic $D$-ball with $p$ internal vertices and its associated $(D+1)$-ary tree, denoted by $M_p$ and $T_p$, respectively. As noted earlier, the ordered set of internal vertices $m_p$ is in bijective correspondence with the ordered set of elementary tree vertices~$t_p$.

 Both of these sets have a natural metric structure inherited from their respective graphs, namely the graph distance:
\begin{Definition}[graph distance] For any connected graph, the \emph{graph distance} $d(v_s,v_t)$ between two vertices~$v_s$ and~$v_t$ is the minimal number edges in any contiguous path journeying from~$v_s$ to~$v_t$.
\end{Definition}
The graph distance is a discrete metric and in order to examine the large-$p$ limit as a continuum limit, we need to extend it to a continuum metric. To do so, one utilizes the order on the set~$m_p$~($t_p$) to arrange pairs $m_p \times m_p$ ($t_p\times t_p$) as the integer points on $[0,p-1]\times[0,p-1]$. One then interpolates between these integer points using a piecewise linear interpolation on the triangles with integer co-ordinates $(v_1,v_2)$, $(v_1+1,v_2)$, $(v_1,v_2+1)$ and $(v_1+1,v_2+1)$, $(v_1+1,v_2)$, $(v_1,v_2+1)$. We denote this metrics by~$d_{m_p}$~($d_{t_p}$).

To retain a compact space in the large-$p$ limit, rather than just letting the structure get inf\/initely large, one requires an $p$-dependent rescaling of the metric. This motivates the factor $\sqrt{(D+1)p/D}$. Of course, it is the exponent of this factor that determines the Hausdorf\/f dimension. Afterall, the volume of the melonic $D$-ball, with respect to~$d_{m_p}$, is proportional to~$p$ and thus, it is the precisely the rescaling $d_{m_p}\rightarrow d_{m_p}/({\rm volume})^{1/d_H}$ required to keep the limiting metric space compact that determines the Hausdorf\/f dimension.

The f\/inal factor $\Lambda_{\Delta}$ requires a little more explanation.

{\bf Calculating $\boldsymbol{\Lambda_{\Delta}}$.} Consider an associated pair $M_p$ and $T_p$, along with their vertex sets~$m_p$ and~$t_p$. Since the vertex sets are in bijective correspondence, we shall refer to their vertices using the same label~$v$. Associated to a rooted graph is the notion of depth:
\begin{Definition}[depth] The \emph{tree depth} of $v$ is the distance, with respect to $d_{t_p}$, of the associated elementary vertex to the root vertex. The \emph{depth} of~$v$ is the distance, with respect to~$d_{m_p}$, of the associated internal vertex to the internal vertex associated to the f\/irst elementary $D$-ball.
\end{Definition}
The tree depth of $v$ is just length of the branch joining it to the root, that is, the number of characters in its associated word~$w$. Denoting the quantity by $\widetilde\Lambda$, one f\/inds for example:
\begin{gather*}
 \widetilde\Lambda(0;10132120312) = 12.
\end{gather*}
The depth, on the other hand, is not so simple to calculate. However, it can also be deciphered from the associated word in the following manner. As explained earlier, a branch within the $(D+1)$-ary tree corresponds to a nested sequence of elementary $D$-balls within the melonic $D$-ball.
As a result of this nesting, the path of minimal graph distance leading from this internal vertex to the initial internal vertex lies within this sequence of elementary $D$-balls. However, the connectivity of the graph causes its value to deviate form that of the tree depth; the insertion of a~$D$-ball does not mean that the new vertex is necessarily further away from the root vertex. In fact, consider a nested sequence of elementary $D$-balls and their associated internal vertices. Pick out the last elementary $D$-ball inserted. Say it has color $i$ and assume that its associated internal vertex is the f\/irst internal vertex in the sequence that is at depth~$r$. All other internal vertices in that sequence have depth less than~$l$.\footnote{Internal vertices inserted earlier certainly do not have depth greater than $r$.} Now, within this latest $D$-ball, insert an arbitrary nested sequence containing at least one elementary $D$-ball of each color \emph{except} some color $j \neq i$. It emerges that the internal vertices associated to this newly inserted sequence all have depth~$r$. Now insert a elementary $D$-ball of color~$j$. Its associated internal vertex has depth~$r+1$. We do not elaborate further on the reasoning here. However, the interested reader can convince oneself quite readily by drawing out the example given here or by looking in~\cite{melondimension}.

Now its a matter of laying out an algorithm for calculating the depth from the associated word. In this regard, let us denote by $W_{D+1}$ the set of words containing every letter of the alphabet $\Sigma_{D+1}=\{0,1,\dots, D\}$ at least once. Consider a vertex $v$ labelled by the word $w=(0 ; u_1u_2\dots u_n)$. The depth of $v$ corresponds to division of $w$ into disjoint adjacent subwords $\tau_r$, comprised of letters of depth~$r$.
Thus, $\tau_0=0$. Then, $\tau_1 = u_1\dots u_{a_1} $, with $ u_1,u_2,\dots, u_{a_1} \neq 0$ and $u_{a_1+1}=0$. Furthermore, $\tau_r$, for $r> 1$, may be one of two forms: \textit{i}) $\tau_r = u_{a_{r-1}+1} \dots u_{a_r} $ such that $\tau_r \notin W_{D+1}$ but $\tau_r u_{a_{r}+1 } \in W_{D+1}$;
\textit{ii}) $\tau_r = u_{a_{r-1}+1} \dots u_n $ if $u_{a_{r-1}+1} \dots u_n \notin W_{D+1}$. This second possibility accounts for the fact that the last subword might be incomplete.

The depth of a vertex $v$ with word $w = \tau_0\tau_1 \dots \tau_k$ is
\begin{gather*}
 \Lambda(w) =k .
\end{gather*}
As an example, take the branch of a $(3+1)$-ary tree illustrated in Fig.~\ref{fig:words}. The division in subwords goes as follows
\begin{gather*}
 w = (0;10132120312) = (0)(1)(013)(2120)(312).
\end{gather*}
Thus, $\Lambda(0;10132120312) = 4$.

The distance, with respect to $d_{m_p}$, between any two vertices can be well estimated from their respective depths with respect to their latest common ancestor. Consider two vertices~$v_1$ and~$v_2$ with associated words~$ww_1$ and~$ww_2$. Thus, the two words have~$w$ in common, but lie on dif\/ferent sub-branches thereafter. The following inequality holds:
\begin{gather*}
 | d_{m_p}(v_1,v_2) - \Lambda(w_1) - \Lambda(w_2) | \le 6.
\end{gather*}
Since the right hand side is $p$-independent, the rescaled inequality becomes increasingly strict as $p$ increases.

Finally for the factor $\Lambda_{\Delta}$, which is just the average ratio $\Lambda(w)/\widetilde\Lambda(w)$ as $p$ becomes large:
\begin{Lemma}\label{lem:rescaling}
Let $u_1,\dots, u_p$ be a sequence of random variables uniformly drawn from $\Sigma_{D+1}$, and denote $w=0u_1\dots u_p$. One has
 \begin{gather*}
 \frac{1}{p} \Lambda(w) \underset{p\rightarrow \infty}{\rightarrow} \Lambda_{\Delta} , \qquad \Lambda_{\Delta}^{-1}
 = (D+1) \sum_{ 0 \le r \le D} (-1)^{D-r} \binom{D}{r} \frac{r}{(D+1-r)^{2}} .
 \end{gather*}
\end{Lemma}
Lemma \ref{lem:rescaling} declares that, on average, the depth of an internal vertex in a melonic $D$-ball is, up to a constant rescaling by~$\Lambda_{\Delta}$, just the tree depth in the associated $(D+1)$-ary tree. Loosely speaking, this factor occurs in Theorem~\ref{th:GH} to take into account the connectivity of the melonic $D$-balls, given that the limiting space has a tree-like structure.

{\bf Melonic $D$-balls as random variables.} In the context of Theorem \ref{th:GH}, $M_p$ is to be viewed as a random variable with uniform distribution upon the space of melonic $D$-balls with $p$ internal vertices. Thus, by association, the metric space $\big(m_p, d_{m_p}/\Lambda_\Delta \sqrt{(D+1)p/D}\big)$ is a~random variable. Thus, the sequence above is a sequence of random variables and convergence means convergence in distribution, that is, the (cumulative) distribution functions associated to the melonic $D$-ball metric spaces random variables converge to that of the continuum random tree in the large-$p$ limit.

Note that the generating function mentioned earlier (as well as the tensor model generating function) weights all melonic $D$-balls with the same number of internal vertices equally. Thus, a uniform distribution is motivated from that context.

{\bf Continuum random tree.} A continuum random tree (CRT) $(\mathcal{T}_{2e}, d_{2e})$ is def\/ined as a rooted real tree encoded by twice
a normalized Brownian excursion $e$ and endowed with a metric $d_{2e}$.

\begin{figure}[t]\centering
\includegraphics[scale =1.2]{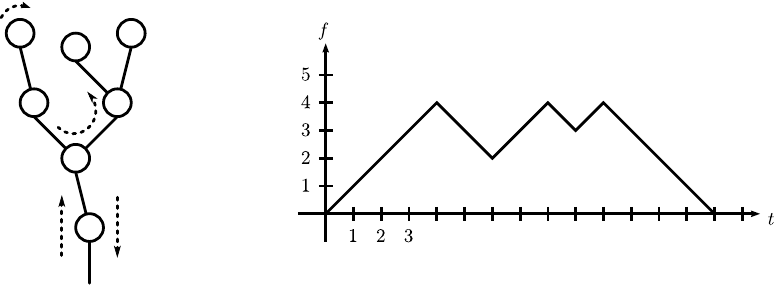}
\caption{A defoliated $(D+1)$-ary tree and its associated contour walk.}\label{fig:excursion}
\end{figure}

One is probably more familiar with rooted discrete trees, of which the colored rooted $(D+1)$-ary trees are examples. Colored rooted $(D+1)$-ary trees (like all discrete trees) have an associated contour walk. Consider such a tree with~$p$ (elementary) vertices. (For simplicity, we shall consider its defoliated version, that is, all leaves removed.) Starting from the base of the tree, one traverses the perimeter of the tree, passing from one vertex to the next in unit time-steps. One considers the following continuous function~$f(t)$, with $f(0) = 0$. As one travels, $f(i) = d_{t_p}(v) + 1$, where~$v$ is the vertex one encounters at the $i$th time-step. (For the value at intermediate times, one linearly interpolates between the time-steps.) The procedure is illustrated in Fig.~\ref{fig:excursion}. Given the construction, one has that the journey ends at time-step $2p$, with $f(2p) = 0$ and $f(t)>0$ for $0< t< 2p$. One has thus associated to any tree some (f\/ixed) walk~$f$. For random trees with~$2p$ vertices, the contour walk becomes a~random walk with~$2p$ steps.

Any real continuous function $f(t)$, such that $f(0) = f(1) = 0$ and $f(t)> 0$ for $0<t<1$, encodes a rooted real tree $\mathcal{T}_f$. To get to the tree, one must set up the following equivalence. For all $s,t\in[0,1]$, set $m_f(s,t) = \inf_{\min(s,t)\leq r \leq \max(s,t)} f(r)$. Then
\begin{gather*}
s\underset{f}{\sim} t \quad\iff\quad f(s) = f(t) = m_f(s,t) .
\end{gather*}
Then, the rooted real tree is the quotient: $\mathcal{T}_f = [0,1]/\underset{f}{\sim}$. The distance on the tree is given by\footnote{In the interests of clarity, let us calculate the distance between two points in the tree of Fig.~\ref{fig:excursion}. Examining the associated contour walk, the heights of the 4th and 9th points encountered are $f(4) = 4$ and $f(9) = 3$, respectively. Meanwhile, the minimal height along the intervening contour is $m_f(4,9) = 2$. Then, the distance between these two points is really the sum of their respective vertical heights above this minimum $d_f(4,9) = (f(4) - m_f(4,9)) + (f(9) - m_f(4,9)) = 3$.}
\begin{gather*}
d_f(s,t) = f(s) + f(t) - 2m_f(s,t) .
\end{gather*}
One can pick out the branching vertices of the tree as those values in $[0,1]$ that are congruent to two or more other values. This real tree dif\/fers from a discrete tree in that one has precise distance information along the edges of the tree.

The Wiener process is a stochastic process $W_t$ (that is a random variable for every time $t$) such that $W_0=0$, $t\to W_t$ is almost surely continuous, $W_t$ has independent increments and $W_t-W_s$ is distributed on a normal distribution of mean $0$ and variance $\sigma^2=t-s$ for $s\le t$. The normalized Brownian excursion $e_t$ is a Wiener process conditioned to be positive for $0<t<1$ and be at $0$ at time $1$. It is formally represented by the path integral measure
\begin{gather*}
 d\mu_{e} = \frac{1}{Z} \bigl[ dq(t)\bigr]\Bigg{|}_{\substack{q(0)=q(1)=0\\ q(t)> 0}} e^{-\frac{1}{2}\int_0^1 [\dot q(t)]^2 dt} ,
\end{gather*}
with $Z$ a normalization constant.

The CRT $(\mathcal{T}_{2e}, d_{2e})$ is the random tree associated to twice a normalized Brownian excur\-sion~$2e$.

\textbf{Gromov--Hausdorf\/f topology and convergence.} Since one considers the convergence of a~sequence of random metric spaces, one must endow the space of metric spaces along with an appropriate topology. This is provided by the Gromov--Hausdorf\/f topology on the space of isometry classes of compact metric spaces.

To begin, one considers a metric space $(E, d_E)$. The Hausdorf\/f distance between two compact sets, $K_1$ and $K_2$, in $E$ is
\begin{gather*}
d_{\textrm{Haus}(E)}(K_1,K_2) = \inf\big\{r \,|\, K_1\subset K_2^r, \, K_2 \subset K_1^r\big\} ,
\end{gather*}
where $K_i^r = \bigcup_{x\in K_i} B_E(x,r)$ is the union of open balls of radius~$r$ centered on the points of~$K_i$.

Now, given two compact metric spaces $( E_i ,d_i)$, the Gromov--Hausdorf\/f distance between them is
\begin{gather*}
d_{\textrm{GH}}(E_1,E_2) = \inf\{d_{\textrm{Haus}(E)}(\phi_1(E_1),\phi_2(E_2)) \} ,
\end{gather*}
where the inf\/imum is taken on all metric spaces $E$ and all isometric embeddings $\phi_1$ and $\phi_2$ from $(E_1,d_1)$ and $(E_2,d_2)$ into $(E,d_E)$.

It emerges that $\mathbb{K}$, the set of all isometry classes of compact metric spaces, endowed with the Gromov--Hausdorf\/f distance $d_{\rm GH}$ is a~complete metric space in its own right. Therefore, one may study the convergence (in distribution) of $\mathbb{K}$-valued random variables.

Of course, the Gromov--Haudorf\/f topology is not the exclusive topology for these metric spaces, but fortuitously, it is well-adapted to the study of quantities that are dependent on the size of the melonic $D$-balls, quantities such as the diameter, the depth, the distance between two random points and so forth.

\subsection{Spectral dimension} \label{sec:spectral}

The third and f\/inal exponent that we would like to examine for melonic $D$-balls is the spectral dimension. The spectral dimension of a manifold is the dimension experienced by a dif\/fusion process and is extracted by evaluating the logarithmic derivative of the return probability with respect to dif\/fusion time
\begin{gather*}
d_S = -2\frac{d\log P(t)}{d \log t}.
\end{gather*}
To make this more concrete, let us consider the simple case of dif\/fusion on $D$-dimensional Euclidean space. The associated heat equation\footnote{In more generality, consider a Riemannian manifold $M$ and a temperature distribution $f_0\colon M \rightarrow \mathbb{R}$. The \emph{heat equation} evolves $f_0$, in the sense that there exists $f\colon M\times (0,T) \rightarrow \mathbb{R}$ that solves
$ \left(\frac{\partial}{\partial t} - \nabla \right) f = 0$,
and $f(x, 0) = f_0(x)$. The \emph{heat kernel} is $K \colon M\times M\times (0,T)\rightarrow \mathbb{R}$, which satisf\/ies
$ \left(\frac{\partial}{\partial t} - \nabla_x \right) K = 0$,
which initial condition
$\lim\limits_{t\rightarrow 0} K(x,y,t) = \delta(x,y)$.}
 gives rise to a \emph{heat kernel}
\begin{gather*}
 K(t,x,y) = \frac{1}{4\pi t^{D/2}}e^{-|x-y|^2/4t}.
\end{gather*}
The spatial coincidence limit gives rise to the \emph{return probability}
\begin{gather*}
 P(t) = \frac{1}{4\pi t^{D/2}},
\end{gather*}
from which the result $d_S = D$ follows readily. Obviously, the spectral dimension coincides with the Hausdorf\/f dimension in this elementary case, but this is not true in general.

{\bf Dif\/fusion on graphs.} On graphical structures, the dif\/fusion process may be modelled using a random walker taking one step per unit time, where a step takes the walker from a vertex to one of its neighbours. The return probability~$P(t)$ is then simply the probability that the walker returns to her starting point at time~$t$.

As one might imagine, the spectral dimension depends non-trivially on the connectivity of the graph. Moreover, the connectivity of the graph can make the computation of the return probability very involved. In our case, we wish to compute the return probability for a generic melonic $D$-ball, which seems rather intractable. So we are motivated to examine the other representations. The spectral dimension of colored rooted $(D+1)$-ary trees is a simple extension of the result in \cite{jw}. However, the rooted melonic graphs bear a closer relationship to the melonic $D$-balls, as they are their topological dual spaces.

\begin{figure}[t]\centering
 \includegraphics[scale=1.2]{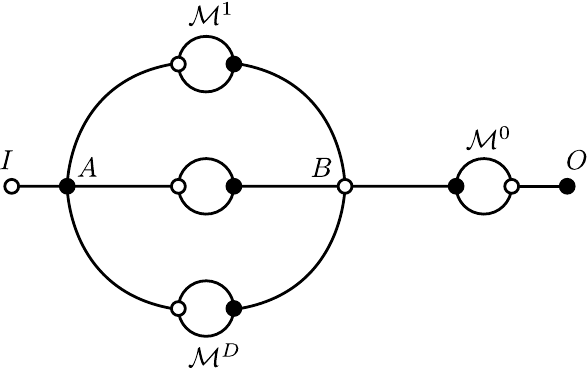}
 \caption{A rooted melonic graph $\mathcal{M}$ with sub-melons $\mathcal{M}^i$.}\label{fig:melongraph}
\end{figure}

To proceed, consider a rooted melonic graph, with external color 0. Such a graph is drawn in Fig.~\ref{fig:melongraph}.

One denotes it by $\mathcal{M}$. It is worth recounting how $\mathcal{M}$ was constructed. One started from an initial elementary melon of color~$0$, that is, two vertices sharing $D$ edges. These two original vertices are denoted by~$A$ and~$B$, Fig.~\ref{fig:melongraph}. Thereafter, one repeatedly inserted elementary melons of various colors. Some were inserted along the active edge of color $0$ at $B$ or nested within earlier insertions. Meanwhile, others were inserted along the active edges of color $i\in\{1,\dots, D\}$ joining~$A$ to~$B$ or nested within earlier insertions. Due to this iterative construction, excising vertices~$A$ and~$B$ results in $D+1$ disconnected subgraphs that are themselves rooted melonic graphs, each with a distinct external color. One denotes this property by ${\mathcal{M}}={\mathcal{M}}^1\cup {\mathcal{M}}^2\cup \cdots\cup {\mathcal{M}}^{D} \cup {\mathcal{M}}^0$, where ${\mathcal{M}}^i$ labelled the rooted melonic graph with external edges of color~$i$.

Any rooted melonic graph has two external vertices, labelled $I$ and $O$. The rest are the familiar~$(D+1)$-valent internal vertices. To calculate the return probability for a generic graph, one should average over the return probabilities attached to each vertex. However, rooted melonic graphs possess a recurrent connectivity structure that renders this unnecessary and ultimately, we wish to investigate return probabilities to the vertex~$I$. To this end, we shall deal with the return probabilities ($I\rightarrow I$ and $O\rightarrow O$) and transit probabilities~($I\rightarrow O$ and $O\rightarrow I$). Utilizing this matrix of probabilities, the recurrent connectivity permits us to set up a~recurrence relation satisf\/ied by these \emph{return/transit probabilities}. We shall detail this presently.

{\bf Generating functions.} As usual, we are interested in properties \emph{on average} (this time over the set of rooted melonic graphs) at large~$p$. Thus, we utilize the convenience of a generating function for the return probabilities to $I$ at a given $t$, over the set of rooted melonic graphs. Moreover, is it simpler to extract the spectral dimension by considering the generating function for return probabilities over the set of all $t$. Thus, we def\/ine the quantity
 \begin{gather*}
 P^{II}(z,y) := \sum_{p = 0}^{\infty}\sum_{t = 0}^{\infty} \widetilde{P}^{II}_{p}(t)y^tz^p ,\qquad\textrm{where}\quad \widetilde{P}^{II}_p(t) := \sum_{\mathcal{M}\colon p_\mathcal{M} = p} P^{II}_{\mathcal{M}}(t).
 \end{gather*}
 We may extract the spectral dimension by examining $P^{II}(z,y)$ in the vicinity of $z \sim z_c$ and $y\sim 1$.
 \begin{Theorem} \label{th:spectral} The return probability generating function on the set of rooted melonic graphs satisfies
 \begin{gather*}
 \frac{\partial P^{II}}{\partial z}(z,y) = (1 - z/z_c)^{\Delta(d_S/2 - 1) - \gamma}
 \left[\Phi\left(\frac{1-y}{(1 - z/z_c)^{\Delta}}\right) + \Phi\left(\frac{1+y}{(1 - z/z_c)^{\Delta}}\right)\right] ,
 \end{gather*}
 within the disc $|y|< 1$, $z_c = D^D/(D+1)^{D+1}$, $\gamma = 1/2$, $\Delta = 3/2$ and
 \begin{gather*}d_S = 4/3.\end{gather*}
 \end{Theorem}

{\bf Random walks on a melon.} Consider a random walk on a rooted melonic graph ${\mathcal{M}}$. The walker takes one step per unit time. If the walker is at one of the external points, one can see from Fig.~\ref{fig:melongraph} that it steps with probability one to its unique neighbor. Meanwhile, if the walker is at any of the internal $(D+1)$-valent vertices, it steps with equal probability, $\frac{1}{D+1}$, to any one of its $D+1$ neighbors.

Thus, any given walk with $t$ steps takes place with a probability that is simply the product of the probabilities from each of its $t$ constituent jumps.

\looseness=-1 {\bf Return/transit random walks.} \emph{Return walks} are walks that start and f\/inish at the external vertex. \emph{Transit walks} start at one external vertex and f\/inish at the other external vertex. Such return/transit random walks have a generating function, the \emph{return/transit probability ge\-ne\-ra\-ting function}: $P^{XY}_{\mathcal{M}}(y)$, where $X \in \{I,O\}$, $Y \in \{I,O\}$ and $y \in \mathbb{C}$. This may be expanded as
 \begin{gather*}
 P^{XY}_{\mathcal{M}}(y) = \sum_{t = 0}^\infty y^t P^{XY}_{\mathcal{M}}(t),
 \end{gather*}
where $P^{XY}_{\mathcal{M}}(t)$ is the probability that the walker arrives at external vertex~$Y$ at time~$t$, given that it starts at external vertex $X$ at time~0. This probability is simply a probability weighted sum over relevant return/transit random walks.

{\bf 1st-return/1st-transit random walks.} \emph{$1$st-return walks} are return walks that spend the intervening period internal vertices.
\emph{1st-transit walks} are transit walks that spend the intervening period at internal vertices. Again such 1st-return/1st-transit random walks have a generating function, the \emph{$1$st-return/$1$st-transit probability generating function}: $P^{1,XY}_{\mathcal{M}}(y)$, where $X \in \{I,O\}$, $Y \in \{I,O\}$ and $y \in \mathbb{C}$. This may be expanded as
 \begin{gather*}
 P^{1,XY}_{\mathcal{M}}(y) = \sum_{t = 0}^\infty y^t P^{1,XY}_{\mathcal{M}}(t),
 \end{gather*}
where $P^{1,XY}_{\mathcal{M}}(t)$ is the probability that the walker arrives at external vertex~$Y$ at time~$t$, given that it starts at external vertex $X$ at time 0 \emph{and} spends the intervening time period at internal vertices. This probability is simply a probability weighted sum over relevant 1st-return/1st-transit random walks.

\begin{Remark}[initial condition] The simplest rooted melonic graph consists of two external vertices connected by a line. It will be denoted by $\mathcal{M}_{(0)}$ and its 1st-return/1st-transit probability generating function is
\begin{gather*}
 P^1_{\mathcal{M}_{(0)}}(y) =
 \begin{pmatrix}
 0 & y \\ y & 0
 \end{pmatrix} .
\end{gather*}
\end{Remark}

{\bf Consistency relations.} There are then two important relations:
\begin{itemize}\itemsep=0pt
 \item \emph{Any return/transit walk can be decomposed as a sequence of $1$st-return/$1$st-transit walks.}
As a result, the return/transit probability generating function can be expressed in terms of the 1st-return/1st-transit probability generating function
\begin{gather*}
 P^{XY}_{\mathcal{M}}(y) = \left[\frac{1}{1 - P^{1}_{\mathcal{M}}(y)}\right]^{XY} .
\end{gather*}

\item \emph{Any $1$st-return/$1$st-transit walks on $\mathcal{M}$ can be decomposed as a sequence of $1$st-return/$1$st-transit walks on its sub-melons $\mathcal{M}^i$, where $i\in\{0,\dots,D\}$.} Again, this relation is mirrored in the generating functions
\begin{gather*} 
 P^1_{ \mathcal{M} } = E^{22} P^1_{ {\mathcal{M}}^0 } E^{22} + \bigl( E^{12} y + E^{22} P^1_{ {\mathcal{M}}^0 } E^{11} \bigr)\\
\hphantom{P^1_{ \mathcal{M} } =}{}
\times \frac{ 1 }{ D+1 - \sum\limits_{i=1}^D P^1_{{\mathcal{M}}^i } - E^{11} P^1_{ {\mathcal{M}}^0 } E^{11} }
\bigl( y E^{21} + E^{11} P^1_{ {\mathcal{M}}^0 } E^{22} \bigr), \\
 P^1_{\mathcal{M}_{(0)}} = \begin{pmatrix}
 0 & y \\ y & 0
 \end{pmatrix},
\end{gather*}
where
\begin{gather*}
 E^{ab}_{\alpha \beta} = \delta^a_{\alpha} \delta^b_{\beta},\qquad a,b,\alpha,\beta \in \{1,2\} .
\end{gather*}
\end{itemize}

{\bf Return/transit probability matrix.} Finally, consider the space of rooted melonic graphs. We are interested in computing the return probability of a random walker, starting from some specif\/ied point, with respect to this set of graphs. We shall specify this point as the external vertex~$I$ for every rooted melonic graph. This return probability is just one element of the \emph{return/transit probability matrix}
\begin{gather*}
 P^{XY}(z,y) = \sum_{\mathcal{M}} P^{XY}_{\mathcal{M}}(y)z^{p_\mathcal{M}} .
\end{gather*}
To prove Theorem~\ref{th:spectral}, one must solve the consistency relations above. For details of the argument, we refer the reader to~\cite{melondimension,jw}.

\section{Conclusion}
We f\/inish up with a quick word about the future study of $(D+1)$-colored graphs. Clearly, we have just scratched the surface here. The combinatorial, topological, algebraic and metric properties deserve a much more extensive analysis. This could have many benef\/its: by analyzing a broader set of graphs, one would hopefully escape the branched polymer phase in the continuum limit. (In fact, such work has begun~\cite{double-scaling+,double-scaling}, but branched polymers have so far proved resilient.)

One such space is the so-called Brownian sphere~\cite{bsphere}. A long-standing open question is to pin down precisely its spectral dimension (there is a large amount of evidence that it equals two). Were the Brownian sphere to occur as the limiting metric space of a suf\/f\/iciently amenable set of colored graphs, then one might be able to tackle the problem as we did above for the continuum random tree.

Of course, from a quantum gravity perspective, one would like to escape branched polymers in favour of some continuum metric space that has physically interesting characteristics, as they have done in the causal dynamical triangulations approach to quantum gravity \cite{cdt}. In particular, one would like to recover a macroscopic 4-dimensional universe in some limit. Were this possible, then the colored graph approach would have a signif\/icant advantage in that one has the bubble Lie algebra (or some subalgebra thereof) at one's disposal, from which one could aim to extract the underlying symmetries of the continuum limit.

\pdfbookmark[1]{References}{ref}
\LastPageEnding

\end{document}